\documentclass[12pt]{article}

\usepackage{latexsym,amssymb}
\usepackage{graphicx}

\newtheorem{theorem}{Theorem}[section]

\newtheorem{example}[theorem]{Example}
\newtheorem{definition}[theorem]{Definition}

\newtheorem{proposition}[theorem]{Proposition}

\newtheorem{lemma}[theorem]{Lemma}
\newtheorem{corollary}[theorem]{Corollary}

\newenvironment{proof}{\medskip\noindent{\it Proof.\ }}{\mbox{$\Box$}\medskip}

\begin{document}

\def\eqnsep{50pt}

\def\fS{\mathfrak S}

\def\cA{\mathcal A}
\def\cD{\mathcal D}
\def\cS{\mathcal S}
\def\cT{\mathcal T}

\def\cplus{\ {{\ } \atop {+}}\ }
\def\cminus{\ {{\ } \atop {-}}\ }
\def\csum{\hbox{c}\hspace{-11pt}\sum}

\title{Permutations Which Avoid 1243 and 2143, \\
Continued Fractions, and Chebyshev Polynomials}
   
\author{
Eric S. Egge \\
Department of Mathematics \\
Gettysburg College\\
Gettysburg, PA  17325  USA \\
\\
eggee@member.ams.org \\
\\
Toufik Mansour \\
LaBRI (UMR 5800), Universit\'e Bordeaux 1, 351 cours de la Lib\'eration,\\
33405 Talence Cedex, France\\
\\
toufik@labri.fr
}

\maketitle

\begin{abstract}
Several authors have examined connections between permutations which avoid 132, continued fractions, and Chebyshev polynomials of the second kind.
In this paper we prove analogues of some of these results for permutations which avoid 1243 and 2143.
Using tools developed to prove these analogues, we give enumerations and generating functions for permutations which avoid 1243, 2143, and certain additional patterns.
We also give generating functions for permutations which avoid 1243 and 2143 and contain certain additional patterns exactly once.
In all cases we express these generating functions in terms of Chebyshev polynomials of the second kind.

\medskip

\end{abstract}

\section{Introduction and Notation}

Let $\fS_n$ denote the set of permutations of $\{1, \ldots, n\}$, written in one-line notation, and suppose $\pi \in \fS_n$ and $\sigma \in \fS_k$.
We say a subsequence of $\pi$ has {\it type} $\sigma$ whenever it has all of the same pairwise comparisons as $\sigma$.
For example, the subsequence 2869 of the permutation 214538769 has type 1324.
We say $\pi$ {\it avoids} $\sigma$ whenever $\pi$ contains no subsequence of type $\sigma$.
For example, the permutation 214538769 avoids 312 and 2413, but it has 2586 as a subsequence so it does not avoid 1243.
If $\pi$ avoids $\sigma$ then $\sigma$ is sometimes called a {\it pattern} or a {\it forbidden subsequence} and $\pi$ is sometimes called a {\it restricted permutation} or a {\it pattern-avoiding permutation}.
In this paper we will be interested in permutations which avoid several patterns, so for any set $R$ of permutations we write $\fS_n(R)$ to denote the elements of $\fS_n$ which avoid every pattern in $R$ and we write $\fS(R)$ to denote the set of all permutations (including the empty permutation) which avoid every pattern in $R$.
When $R = \{\pi_1, \pi_2, \ldots, \pi_r\}$ we often write $\fS_n(R) = \fS_n(\pi_1, \pi_2, \dots, \pi_r)$ and $\fS(R) = \fS(\pi_1, \pi_2, \dots, \pi_r)$.

Several authors \cite{BCS,JR,Krattenthaler,MansourVainshtein4,RWZ} have shown that generating functions for $\fS(132)$ with respect to the number of subsequences of type $12\ldots k$, for various collections of values of $k$, can be expressed as continued fractions.
The most general result along these lines, which appears as \cite[Theorem 1]{BCS}, states that
\begin{equation}
\label{eqn:introcf}
\sum_{\pi \in \fS(132)} \prod_{k \ge 1} x_k^{\tau_k(\pi)} = \frac{1}{\displaystyle
1 - \frac{x_1}{\displaystyle
1 - \frac{x_1 x_2}{\displaystyle
1 - \frac{x_1 x_2^2 x_3}{\displaystyle
1 - \frac{x_1 x_2^3 x_3^3 x_4}{\displaystyle
1 - \frac{x_1 x_2^4 x_3^6 x_4^4 x_5}{\displaystyle
1 - \cdots}}}}}}.
\end{equation}
Here $\tau_k(\pi)$ is the number of subsequences of type $12\ldots k$ in $\pi$.
Generating functions for $\fS(132)$ have also been found to be expressible in terms of Chebyshev polynomials of the second kind \cite{CW,Krattenthaler,MansourVainshtein4,MansourVainshtein}.
One result along these lines, which appears as \cite[Theorem 2]{Krattenthaler}, \cite[Theorem 3.1]{MansourVainshtein4}, and \cite[Theorem 3.6, second case]{CW}, states that
\begin{equation}
\label{eqn:introCheby}
\sum_{n=0}^\infty |\fS_n(132, 12\ldots k)| x^n = \frac{U_{k-1}\left( \frac{1}{2 \sqrt{x}}\right)}{\sqrt{x} U_k \left( \frac{1}{2 \sqrt{x}} \right)}.
\end{equation}
Here $U_n(x)$ is the $n$th Chebyshev polynomial of the second kind, which may be defined by ${\displaystyle U_n(\cos t) = \frac{\sin((n+1) t)}{\sin t}}$.
Another result along these lines, which appears as \cite[Theorem 3]{Krattenthaler}, states that
\begin{equation}
\label{eqn:intro12kr}
\sum_{\pi} x^{|\pi|} = \sum \prod_{i = 2}^b {{l_{i-1} + l_i - 1}\choose {l_i}} \frac{\left( U_{k-1}\left(\frac{1}{2 \sqrt{x}}\right)\right)^{l_1 - 1}}{\left(U_k\left( \frac{1}{2\sqrt{x}}\right)\right)^{l_1 + 1}} x^{\frac{1}{2}( l_1 - 1 ) + \sum\limits_{j=2}^b l_j}.
\end{equation}
Here the sum on the left is over all permutations in $\fS(132)$ which contain exactly $r$ subsequences of type $12\ldots k$, the quantity $|\pi|$ is the length of $\pi$, and the sum on the right is over all sequences $l_1, l_2, \ldots, l_b$ of nonnegative integers such that $\sum_{i=1}^b l_i {{k+i-2}\choose{k-1}} = r$.
For other results involving $\fS(132)$ and continued fractions or Chebyshev polynomials, see \cite{MansourVainshtein5} and the references therein.

Permutations which avoid 1243 and 2143 are known to have many properties which are analogous to properties of permutations which avoid 132.
For instance, it is well known that $|\fS_n(132)| = C_n$ for all $n \ge 0$, where $C_n$ is the $n$th Catalan number, which may be defined by $C_0 = 1$ and 
$$C_n = \sum_{i=1}^n C_{i-1} C_{n-i} \hspace{30pt} (n \ge 1).$$
(The Catalan number $C_n$ may also be defined by $C_n = \frac{1}{n+1} {{2n} \choose {n}}$.)
As a result, for all $n \ge 0$, the set $\fS_n(132)$ is in bijection with the set of {\it Catalan paths}.
These are the lattice paths from $(0,0)$ to $(n,n)$ which contain only east $(1,0)$ and north $(0,1)$ steps and which do not pass below the line $y = x$.
Kremer \cite[Corollary 9]{Kremer1} has shown that $|\fS_n(1243, 2143)| = r_{n-1}$ for all $n \ge 1$, where $r_n$ is the $n$th Schr\"oder number, which may be defined by $r_0 = 1$ and
$$r_n = r_{n-1} + \sum_{i=1}^{n} r_{i-1} r_{n-i} \hspace{30pt} (n \ge 1).$$
As a result, for all $n \ge 0$, the set $\fS_{n+1}(1243, 2143)$ is in bijection with the set $\cS_n$ of {\it Schr\"oder paths}.
These are the lattice paths from $(0,0)$ to $(n,n)$ which contain only east $(1,0)$, north $(0,1)$, and diagonal $(1,1)$ steps and which do not pass below the line $y = x$.
We write $\cS$ to denote the set of all Schr\"oder paths (including the empty path).
In view of this relationship, we refer to permutations which avoid 1243 and 2143 as {\it Schr\"oder permutations}.
(For more information on pattern-avoiding permutations counted by the Schr\"oder numbers, see \cite{BDPP,Gire,Kremer1,WestCatalanSchroder}.
For generalizations of some of these results, see \cite{BDPP,MansourThesis}.
For a partial list of other combinatorial objects counted by the Schr\"oder numbers, see \cite[pp.  239--240]{StanleyVol2}.)

Motivated by the parallels between $\fS(132)$ and $\fS(1243, 2143)$, in this paper we prove analogues of (\ref{eqn:introcf}), (\ref{eqn:introCheby}), (\ref{eqn:intro12kr}), and several similar results for $\fS(1243, 2143)$.
We begin with some results concerning $\fS(1243, 2143)$ and continued fractions.
We first define statistics $\tau_k$, $k \ge 1$, on $\cS$ and $\fS(1243, 2143)$.
On $\fS(1243, 2143)$, the statistic $\tau_k$ is simply the number of subsequences of type $12\ldots k$.
On $\cS$, the statistic $\tau_k$ is a sum of binomial coefficients over east and diagonal steps.
We then give a combinatorial definition of a bijection $\varphi : \cS \rightarrow \fS(1243, 2143)$ with the property that $\tau_k(\varphi(\pi)) = \tau_k(\pi)$ for all $k \ge 1$ and all $\pi \in \cS$.
Using $\varphi$ and a result of Flajolet \cite[Theorem 1]{Flajolet}, we prove the following analogue of (\ref{eqn:introcf}).
\begin{equation}
\label{eqn:introscf}
\sum_{\pi \in \fS(1243, 2143)} \prod_{k \ge 1} x_k^{\tau_k(\pi)} = 
1 + \frac{x_1}{\displaystyle 1 - x_1 - \frac{x_1 x_2}{\displaystyle
1 - x_1 x_2 - \frac{x_1 x_2^2 x_3}{\displaystyle
1 - x_1 x_2^2 x_3 - \cdots}}}.
\end{equation}
By specializing the $x_i$s in (\ref{eqn:introscf}), we obtain continued fraction expansions for several other statistics on $\fS(1243, 2143)$.
In particular, we show that for all $k \ge 1$,
$$\sum_{n=0}^\infty |\fS_n(1243, 2143, 12\ldots k)| x^n = 1 + \frac{x}{\displaystyle
1 - x - \frac{x}{\displaystyle
1 - x - \frac{x}{\displaystyle
1 - x \cdots}}}.$$
Here the continued fraction on the right has $k - 1$ denominators.
Following \cite{BCS}, we then define a {\it Schr\"oder continued fraction} to be a continued fraction of the form
$$1 + \frac{m_0}{\displaystyle
1 - m_0 - \frac{m_1}{\displaystyle
1 - m_1 - \frac{m_2}{\displaystyle
1 - m_2 - \frac{m_3}{\displaystyle
1 - m_3 - \cdots}}}},$$
where $m_i$ is a finite monic monomial in a given set of variables for all $i \ge 0$.
We prove that the multivariate generating function for a countable family of statistics on $\fS(1243, 2143)$ can be expressed as a Schr\"oder continued fraction if and only if each statistic is a (possibly infinite) linear combination of the $\tau_k$s and each $\tau_k$ appears in only finitely many of these linear combinations.
This result is an analogue of \cite[Theorem 2]{BCS}.

We then turn our attention to analogues of (\ref{eqn:introCheby}) and (\ref{eqn:intro12kr}).
For any $k \ge 2$ and any $\sigma \in \fS_{k-1}$ we give the generating function for $|\fS_n(1243, 2143, k \sigma)|$ in terms of the generating function for $|\fS_n(1243, 2143, \sigma)|$.
Using this result, we show that
$$\sum_{n=0}^\infty |\fS_n(1243, 2143, 12\ldots k)| x^n = 1 + \frac{\sqrt{x} U_{k-2}\left(\frac{1-x}{2 \sqrt{x}}\right)}{U_{k-1}\left( \frac{1-x}{2 \sqrt{x}}\right)}$$
and
$$\sum_{n=0}^\infty |\fS_n(1243, 2143, 213\ldots k)| x^n = 1 + \frac{\sqrt{x} U_{k-2}\left(\frac{1-x}{2 \sqrt{x}}\right)}{U_{k-1}\left( \frac{1-x}{2 \sqrt{x}}\right)}$$
for all $k \ge 1$.
Both of these results are analogues of (\ref{eqn:introCheby}).
We then use $\varphi$ and some well-known results concerning lattice paths to show that
$$\sum_{\pi} x^{|\pi|} = \sum \prod_{i=0}^b {{l_i + l_{i+1} + m_i - 1} \choose {l_{i+1} + m_i}} {{l_{i+1} + m_i} \choose {m_i}} \frac{\left(U_{k-2}\left(\frac{1-x}{2\sqrt{x}}\right)\right)^{l_0 - 1}}{\left( U_{k-1}\left(\frac{1-x}{2\sqrt{x}}\right)\right)^{l_0+1}} x^{\frac{1}{2}(1-l_0) + \sum\limits_{j=0}^b ( l_j + m_j )}.$$
Here the sum on the left is over all permutations in $\fS(1243, 2143)$ which contain exactly $r$ subsequences of type $12 \ldots k$ and the sum on the right is over all sequences $l_0, l_1, \ldots, l_b,$ and $m_0, m_1, \ldots, m_b$ of nonnegative integers such that $r = \sum_{i=0}^b (l_i + m_i) {{k+i-1} \choose {k-1}}.$
This result is an analogue of (\ref{eqn:intro12kr}).

In the next two sections of the paper we give enumerations and generating functions for various sets of permutations in $\fS(1243, 2143)$.
For instance, we show that
\begin{equation}
\label{eqn:intro231gf}
\sum_{n=0}^\infty |\fS_n(1243, 2143, 2134\ldots k)| x^n = 1 + x \frac{f_{k-1}(x)}{f_k(x)}
\end{equation}
and
\begin{equation}
\label{eqn:intro321gf}
\sum_{n=0}^\infty |\fS_n(1243, 2143, 3214\ldots k)| x^n = 1 + x \frac{g_{k-1}(x)}{g_k(x)}
\end{equation}
for all $k \ge 3$.
Here $f_2(x) = (x-1)^2$,
$$f_k(x) = (1-2x)^2 (\sqrt{x})^{k-3} U_{k-3}\left(\frac{1-x}{2\sqrt{x}}\right) - (1-x)^2  (\sqrt{x})^{k-2} U_{k-4}\left(\frac{1-x}{2\sqrt{x}}\right)$$
for all $k \ge 3$, and
$$g_k(x) = -(1+2x-x^2)(\sqrt{x})^{k+2} U_k\left( \frac{1-x}{2\sqrt{x}}\right) + (x^4-4x^3+2x^2+1)(\sqrt{x})^{k-1} U_{k-1}\left(\frac{1-x}{2\sqrt{x}}\right)$$
for all $k \ge 2$.
Setting $k = 3$ in (\ref{eqn:intro231gf}) and (\ref{eqn:intro321gf}), we find that
\begin{equation}
\label{eqn:intro231}
|\fS_n(1243, 2143, 231)| = (n+2) 2^{n-3} \hspace{30pt} (n \ge 2)
\end{equation}
and
\begin{equation}
\label{eqn:intro321}
|\fS_n(1243, 2143, 321)| = {{n-1}\choose{0}} + {{n-1}\choose{1}} + 2 {{n-1}\choose{2}} + 2 {{n-1}\choose{3}} \hspace{30pt} (n \ge 1).
\end{equation}
It is an open problem to provide combinatorial proofs of (\ref{eqn:intro231}) and (\ref{eqn:intro321}).
We also show that 
$$\sum_{\pi} x^{|\pi|} = \frac{x (1+x) (1-x)^2}{\left(U_{k-1}\left(\frac{1-x}{2\sqrt{x}}\right)\right)^2},$$
where the sum on the left is over all permutations in $\fS(1243, 2143)$ which contain exactly one subsequence of type $213\ldots k$.
It is an open problem to give the sum $\sum\limits_{\pi} x^{|\pi|}$ in closed form when it is over all permutations in $\fS(1243, 2143)$ which contain exactly $r$ subsequences of type $213\ldots k$, where $r \ge 2$.

We conclude the paper by collecting several open problems related to this work.

\section{Statistics and a Product for Schr\"oder Paths}

In this section we define a family of statistics on Schr\"oder paths.
We then recall the first-return product on Schr\"oder paths and describe the behavior of our statistics with respect to this product.
We begin by recalling the height of an east or diagonal step in a Schr\"oder path.

\begin{definition}
Let $\pi$ denote a Schr\"oder path, let $s$ denote a step in $\pi$ which is either east or diagonal, and let $(x,y)$ denote the coordinates of the left-most point in $\pi$.
We define the {\em height} of $s$, written $ht(s)$, by setting $ht(s) = y - x$.
\end{definition}

We now define our family of statistics on Schr\"oder paths.

\begin{definition}
For any Schr\"oder path $\pi$ and any positive integer $k$ we write
\begin{equation}
\label{eqn:skpath}
\tau_k(\pi) = {{0} \choose {k-1}} + \sum_{s \in \pi} {{ht(s)} \choose {k-1}},
\end{equation}
where the sum on the right is over all east and diagonal steps in $\pi$.
Here we take ${{i} \choose {j}} = 0$ whenever $j < 0$ or $i < j$.
For notational convenience we set $\tau_0(\pi) = 0$ for any Schr\"oder path $\pi$.
\end{definition}

\begin{example}
\label{exmp:schroderpath}
Let $\pi$ denote the Schr\"oder path given in Figure 1, so that $\pi$ is given by $\pi = NDEDNNNNDNEENEDEEE$.
Then $\tau_1(\pi) = 12$, $\tau_2(\pi) = 28$, $\tau_3(\pi) = 35$, $\tau_4(\pi) = 24$, $\tau_5(\pi) = 8$, $\tau_6(\pi) = 1$, and $\tau_k(\pi) = 0$ for all $k \ge 7$.
\end{example}

\begin{figure}[ht]
\begin{center}
\includegraphics[scale=.35]{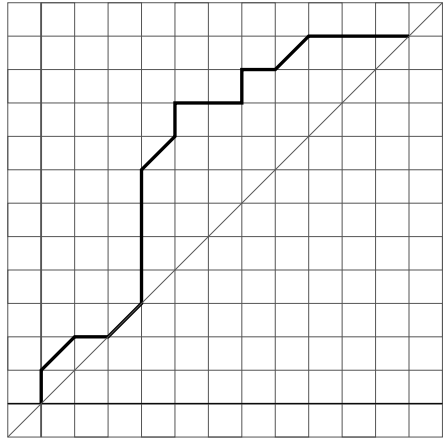}
\caption{The Schr\"oder path of Example \ref{exmp:schroderpath}.}
\end{center}
\end{figure}

Before we recall the first-return product for Schr\"oder paths, we make an observation regarding those paths in $\cS_n$ which begin with a diagonal step.

\begin{proposition}
\label{prop:pathnrecurrence}
\renewcommand\labelenumi{{\upshape (\roman{enumi}) }}
\begin{enumerate}
\item
For all $n \ge 1$, the map
$$
\begin{array}{ccc}
\cS_{n-1} &\longrightarrow& \cS_n \\
\pi &\mapsto& D, \pi \\
\end{array}
$$
is a bijection between $\cS_{n-1}$ and the set of Schr\"oder paths in $\cS_n$ which begin with a diagonal step.
\item $\tau_1(D,\pi) = 1 + \tau_1(\pi)$ for all $\pi \in \cS$.
\item $\tau_k(D,\pi) = \tau_k(\pi)$ for all $k \ge 2$ and all $\pi \in \cS$.
\end{enumerate}
\end{proposition}
\begin{proof}
(i)
This is immediate.

(ii),(iii)
From (\ref{eqn:skpath}) we find that for all $k \ge 1$,
\begin{eqnarray*}
\tau_k(D,\pi) &=& {{0} \choose {k-1}} + {{0}\choose{k-1}} + \sum_{s \in \pi} {{ht(s)} \choose {k-1}} \\[2ex]
&=& {{0} \choose {k-1}} + \tau_k(\pi).
\end{eqnarray*}
Now (ii) and (iii) follow.
\end{proof}

We now define the first-return product on Schr\"oder paths.

\begin{definition}
For any Schr\"oder paths $\pi_1$ and $\pi_2$ we write
$$\pi_1 * \pi_2 = N \pi_1 E \pi_2.$$
\end{definition}

\begin{proposition}
\label{prop:pathinirecurrence}
Let $i$ and $n$ denote positive integers such that $1 \le i \le n$.
Then the following hold.
\renewcommand\labelenumi{{\upshape (\roman{enumi}) }}
\begin{enumerate}
\item The map
$$
\begin{array}{ccc}
\cS_{i-1} \times \cS_{n-i} &\longrightarrow& \cS_n \\
(\pi_1, \pi_2) &\mapsto& \pi_1 * \pi_2 \\
\end{array}
$$
is a bijection between $\cS_{i-1} \times \cS_{n-i}$ and the set of Schr\"oder paths in $\cS_n$ which begin with a north step and first touch the line $y = x$ at $(i,i)$.
\item
For all $k \ge 1$, all $\pi_1 \in \cS_{i-1}$, and all $\pi_2 \in \cS_{n-i}$ we have
\begin{equation}
\label{eqn:taupath}
\tau_k(\pi_1 * \pi_2) = \tau_k(\pi_1) + \tau_{k-1}(\pi_1) + \tau_k(\pi_2).
\end{equation}
\end{enumerate}
\end{proposition}
\begin{proof}
(i)
This is immediate.

(ii)
Fix $k \ge 1$.
Using (\ref{eqn:skpath}) we have
\begin{eqnarray*}
\lefteqn{\tau_k(\pi_1 * \pi_2) = {{0} \choose {k-1}} + \sum_{s \in \pi_1} {{ht(s) + 1} \choose {k-1}} + {{1}\choose{k-1}} + \sum_{s \in \pi_2} {{ht(s)}\choose{k-1}}} & & \\[2.5ex]
&=& {{0}\choose{k-1}} + \sum_{s \in \pi_1} {{ht(s)}\choose{k-1}} + \sum_{s \in \pi_1} {{ht(s)}\choose{k-2}} + {{0} \choose {k-2}} + {{0} \choose {k-1}} + \sum_{s \in \pi_2} {{ht(s)}\choose{k-1}} \\[2.5ex]
&=& \tau_k(\pi_1) + \tau_{k-1}(\pi_1) + \tau_k(\pi_2),
\end{eqnarray*}
as desired.
\end{proof}

\section{Statistics and a Product for Schr\"oder Permutations}
\label{sec:permproducts}

In this section we define a natural family of statistics on Schr\"oder permutations which is analogous to the family of statistics we have defined on Schr\"oder paths.
We then describe a ``product'' on Schr\"oder permutations which behaves nicely with respect to our statistics.
This product is analogous to the first-return product for Schr\"oder paths given in the previous section.
We begin with our family of statistics.

\begin{definition}
For any positive integer $k$ and any permutation $\pi$, we write $\tau_k(\pi)$ to denote the number of increasing subsequences of length $k$ which are contained in $\pi$.
For notational convenience we set $\tau_0(\pi) = 0$ for any permutation $\pi$.
\end{definition}

Observe that for any permutation $\pi$, the quantity $\tau_1(\pi)$ is the length of $\pi$;  we sometimes write $|\pi|$ to denote this quantity.

\begin{example}
If $\pi = 71824356$ then $\tau_1(\pi) = 8$, $\tau_2(\pi) = 16$, $\tau_3(\pi) = 16$, $\tau_4(\pi) = 9$, $\tau_5(\pi) = 2$, and $\tau_k(\pi) = 0$ for all $k \ge 6$.
\end{example}

Observe that we have now defined $\tau_k(\pi)$ when $\pi$ is a Schr\"oder permutation and when $\pi$ is a Schr\"oder path.
This will not cause confusion, however, since it will always be clear from the context which definition is intended.

Before we describe our product for Schr\"oder permutations, we make an observation regarding those Schr\"oder permutations whose largest element appears first.

\begin{proposition}
\label{prop:permnrecurrence}
\renewcommand\labelenumi{{\upshape (\roman{enumi}) }}
\begin{enumerate}
\item
For all $n \ge 1$, the map
$$
\begin{array}{ccc}
\fS_{n-1}(1243, 2143) &\longrightarrow& \fS_n(1243, 2143) \\
\pi &\mapsto& n,\pi \\
\end{array}
$$
is a bijection between $\fS_{n-1}(1243, 2143)$ and the set of permutations in $\fS_n(1243, 2143)$ which begin with $n$.
\item $\tau_1(n,\pi) = 1 + \tau_1(\pi)$ for all $n \ge 1$ and all $\pi \in \fS_{n-1}(1243, 2143)$.
\item $\tau_k(n,\pi) = \tau_k(\pi)$ for all $k \ge 2$, all $n \ge 1$, and all $\pi \in \fS_{n-1}(1243, 2143)$.
\end{enumerate}
\end{proposition}
\begin{proof}
(i)
It is clear that the given map is one-to-one and that if $n,\pi$ is a permutation in $\fS_n(1243, 2143)$ then $\pi \in \fS_{n-1}(1243, 2143)$, so it is sufficient to show that if $\pi \in \fS_{n-1}(1243, 2143)$ then $n,\pi$ avoids $1243$ and $2143$.
To this end, suppose $\pi \in \fS_{n-1}(1243, 2143)$.
Since $\pi$ avoids $1243$ and $2143$, in any pattern of either type in $n,\pi$ the $n$ must play the role of the 4.
But this is impossible, since the $n$ is the left-most element of $n,\pi$, but 4 is not the left-most element of 1243 or 2143.
Therefore $n,\pi$ avoids 1243 and 2143.

(ii)
This is immediate, since $\tau_1(\pi)$ is the length of $\pi$ for any permutation $\pi$.

(iii)
Since $n$ is both the largest and the left-most element in $n,\pi$, it cannot participate in an increasing subsequence of length two or more.
Therefore any such subsequence in $n,\pi$ is contained in $\pi$, and (iii) follows. 
\end{proof}

We now describe our product for Schr\"oder permutations.
To do so, we first set some notation.

Let $\pi_1$ and $\pi_2$ denote nonempty Schr\"oder permutations.
We write $\tilde{\pi_1}$ to denote the sequence obtained by adding $|\pi_2| - 1$ to every entry in $\pi_1$ and then replacing $|\pi_2|$ (the smallest entry in the resulting sequence) with the left-most entry of $\pi_2$.
We observe that $\tilde{\pi_1}$ has type $\pi_1$.
We write $\tilde{\pi_2}$ to denote the sequence obtained from $\pi_2$ by removing its left-most element.

\begin{definition}
\label{defn:permstar}
For any nonempty Schr\"oder permutations $\pi_1$ and $\pi_2$, we write
$$\pi_1 * \pi_2 = \tilde{\pi_1} n \tilde{\pi_2},$$
where $n = |\pi_1| + |\pi_2|$ and $\tilde{\pi_1}$ and $\tilde{\pi_2}$ are the sequences described in the previous paragraph.
\end{definition}

\begin{example}
If $\pi_1 = 3124$ and $\pi_2 = 15342$ then $\tilde{\pi_1} = 7168$, $\tilde{\pi_2} = 5342$, and $\pi_1 * \pi_2 = 716895342$.
\end{example}

\begin{proposition}
\label{prop:perminirecurrence}
Let $i$ and $n$ denote positive integers such that $1 \le i \le n-1$.
Then the following hold.
\renewcommand\labelenumi{{\upshape (\roman{enumi}) }}
\begin{enumerate}
\item
The map
$$
\begin{array}{ccc}
\fS_i(1243, 2143) \times \fS_{n-i}(1243, 2143) &\longrightarrow& \fS_n(1243, 2143) \\
(\pi_1, \pi_2) &\mapsto& \pi_1 * \pi_2 \\
\end{array}
$$
is a bijection between $\fS_i(1243, 2143) \times \fS_{n-i}(1243, 2143)$ and the set of permutations in $\fS_n(1243, 2143)$ for which $\pi(i+1) = n$.
\item 
For all $k \ge 1$, all $\pi_1 \in \fS_i(1243, 2143)$, and all $\pi_2 \in \fS_{n-i}(1243, 2143)$ we have
\begin{equation}
\label{eqn:tauperm}
\tau_k(\pi_1 * \pi_2) = \tau_k(\pi_1) + \tau_{k-1}(\pi_1) + \tau_k(\pi_2).
\end{equation}
\end{enumerate}
\end{proposition}
\begin{proof}
(i)
It is routine to verify that if $\pi_1 \in \fS_i(1243, 2143)$ and $\pi_2 \in \fS_{n-i}(1243, 2143)$ then $\pi_1 * \pi_2 \in \fS_n(1243, 2143)$, so it is sufficient to show that the given map is a bijection.
We do this by describing its inverse.

Suppose $\pi \in \fS_n(1243, 2143)$ and $\pi(i+1) = n$.
Let $f_1(\pi)$ denote the type of the subsequence $\pi(1), \pi(2), \ldots, \pi(i)$ of $\pi$ and let $f_2(\pi)$ denote the permutation of $1,2,\ldots, n-i$ which appears in $\pi$.
Since $\pi \in \fS_n(1243, 2143)$ and $\pi$ contains subsequences of type $f_1(\pi)$ and $f_2(\pi)$ we find that $f_1(\pi) \in \fS_i(1243, 2143)$ and $f_2(\pi) \in \fS_{n-i}(1243, 2143)$.
We now show that the map $\pi \mapsto (f_1(\pi), f_2(\pi))$ is the inverse of the map $(\pi_1, \pi_2) \mapsto \pi_1 * \pi_2$.

It is clear from the construction of $\pi_1 * \pi_2$ that $f_1(\pi_1 * \pi_2) = \pi_1$ and $f_2(\pi_1 * \pi_2) = \pi_2$, so it remains to show that $f_1(\pi) * f_2(\pi) = \pi$.
To this end, observe that since $\pi$ avoids 1243 and 2143 and has $\pi(i+1) = n$, exactly one of $1,2,\ldots, n-i$ appears to the left of $n$ in $\pi$.
(If two or more of $1,2,\ldots,n-i$ appeared to the left of $n$ then two of these elements, together with $n$ and some element to the right of $n$, would form a subsequence of type 1243 or 2143.)
Therefore, the remaining elements among $1,2,\ldots,n-i$ are exactly those elements which appear to the right of $n$.
It follows that $f_1(\pi) * f_2(\pi) = \pi$. 
Therefore the map $\pi \mapsto (f_1(\pi), f_2(\pi))$ is the inverse of the map $(\pi_1, \pi_2) \mapsto \pi_1 * \pi_2$, so the latter is a bijection, as desired.

(ii)
Let $\tilde{\pi_1}$ be as in the paragraph above Definition \ref{defn:permstar}.
Observe that since $\tilde{\pi_1}$ consists of exactly those elements of $\pi_1 * \pi_2$ which are to the left of $n$, there is a one-to-one correspondence between increasing subsequences of length $k-1$ in $\tilde{\pi_1}$ and increasing subsequences of length $k$ in $\pi_1 * \pi_2$ which involve $n$.
Since $\tilde{\pi_1}$ and $\pi_1$ have the same type, there are $\tau_{k-1}(\pi_1)$ of these subsequences.
Now observe that if an increasing subsequence of length $k$ in $\pi_1 * \pi_2$ does not involve $n$, and involves an element of $\tilde{\pi_1}$ other than the smallest element, then it is entirely contained in $\tilde{\pi_1}$.
Similarly, observe that if an increasing subsequence of length $k$ in $\pi_1 * \pi_2$ involves an element of $\pi_2$ other than the left-most element, then it is entirely contained in $\pi_2$.
Therefore every increasing subsequence of length $k$ in $\pi_1 * \pi_2$ which does not involve $n$ is an increasing subsequence of length $k$ in $\tilde{\pi_1}$ or in $\pi_2$.
Since $\tilde{\pi_1}$ and $\pi_1$ have the same type, there are $\tau_k(\pi_1) + \tau_k(\pi_2)$ of these subsequences.
Now (ii) follows.
\end{proof}

Although the results we have given in this section are sufficient for our current purposes, we remark that there are more general results along the same lines.
For example, following \cite{BDPP} and \cite{MansourThesis}, let $T_k$ ($k \ge 3$) denote the set of all permutations in $\fS_k$ which end with $k, k-1$.
Observe that $T_3 = \{132\}$ and $T_4 = \{1243, 2143\}$.
Then there are natural analogues of all of the results in this section for $\fS(T_k)$, where $k \ge 5$.

\section{A Bijection Between $\cS_n$ and $\fS_{n+1}(1243, 2143)$}

Comparing Propositions \ref{prop:permnrecurrence} and \ref{prop:perminirecurrence} with Propositions \ref{prop:pathnrecurrence} and \ref{prop:pathinirecurrence} respectively, we see that for all $n \ge 0$ there exists a bijection $\varphi : \cS_n \rightarrow \fS_{n+1}(1243, 2143)$ such that $\tau_k(\pi) = \tau_k(\varphi(\pi))$ for all $\pi \in \cS_n$.
So far, we have only seen how to compute this bijection recursively.
In this section we use techniques of Bandlow and Killpatrick \cite{BK} and Bandlow, Egge, and Killpatrick \cite{BEK} to compute this bijection directly.

To define our bijection, we first need to introduce some notation.
For all $n \ge 0$ and all $i$ such that $1 \le i \le n-1$, we write $s_i$ to denote the map from $\fS_n$ to $\fS_n$ which acts by interchanging the elements in positions $i$ and $i+1$ of the given permutation.
For example, $s_1(354126) = 534126$ and $s_4(354126) = 354216$.
We apply these maps from right to left, so that $s_i s_j(\pi) = s_i(s_j(\pi))$.

Suppose $\pi \in \cS_n$.
We now describe how to construct the image $\varphi(\pi)$ of $\pi$ under our bijection $\varphi$.
To illustrate the procedure, we give a running example in which the Schr\"oder path $\pi$ is given by $\pi = NDNNEEENNDENEE$, which is illustrated in Figure 2 below.

\begin{figure}[ht]
\begin{center}
\includegraphics[scale=.4]{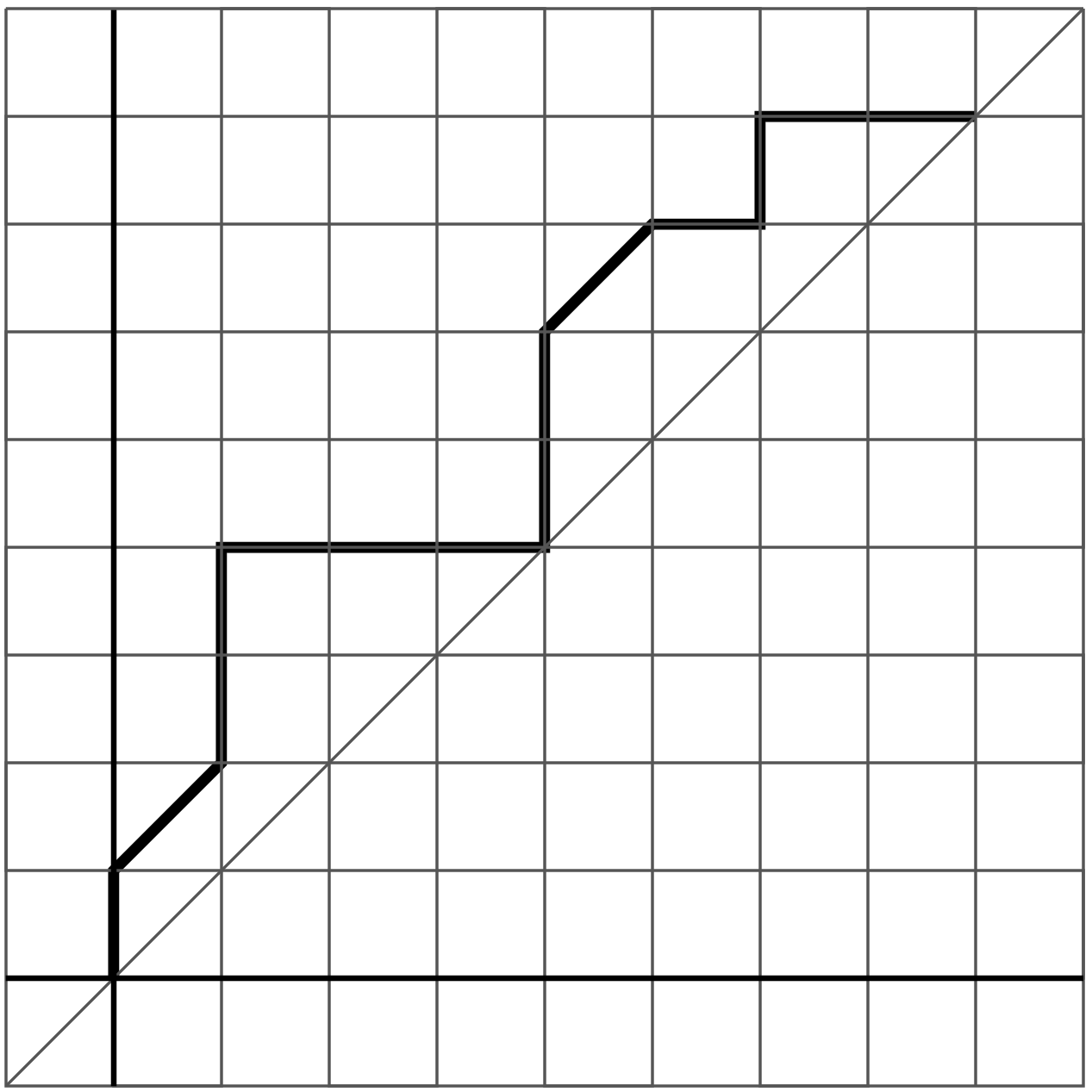}
\caption{The Schr\"oder path $\pi$ used in the running example.}
\end{center}
\end{figure}

To begin, label each upper triangle (i.e. each triangle whose vertices have coordinates of the form $(i-1,j-1)$, $(i-1,j)$, $(i,j)$) which is below $\pi$ and above the line $y = x$ with an $s_i$, where $i$ is the $x$-coordinate of the upper right corner of the triangle.
In Figure 3 below we have labeled the appropriate triangles for $\pi$ with the subscripts for the $s_i$s.
\begin{figure}[ht]
\begin{center}
\includegraphics[scale=.6]{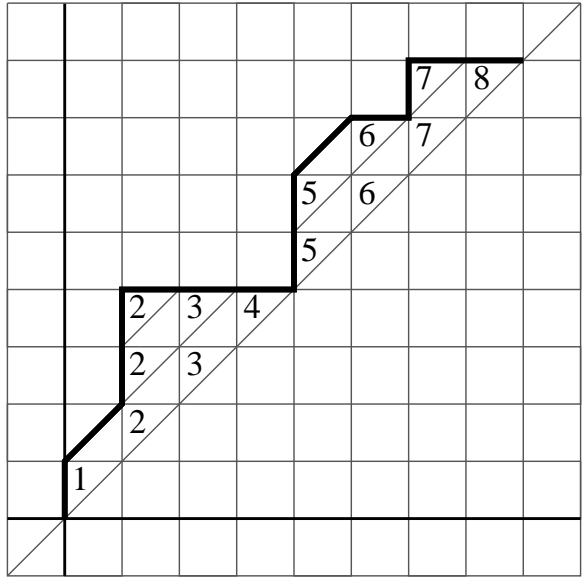}
\caption{The Schr\"oder path $\pi$ with triangles labeled.}
\end{center}
\end{figure}
Now let $\sigma_1$ denote the sequence of $s_i$s which begins with the $s_i$ furthest up and to the right and extends diagonally to the lower left, immediately above the line $y = x$, proceeding until it reaches a north step.
Let $\sigma_2$ denote the sequence of $s_i$s which begins immediately to the left of the beginning of $\sigma_1$ and extends diagonally to the lower left, immediately above $\sigma_1$, until it reaches a north step.
Construct $\sigma_3, \sigma_4, \ldots$ in this fashion until all of the $s_i$s above and to the left of $\sigma_1$ have been read.
Repeat this process with the part of the path below the last east step before the north step which ended $\sigma_1$.
In this way we obtain a sequence $\sigma_1, \ldots, \sigma_k$ of maps on $\fS_n$.
In Figure 4 below we have indicated the sequences $\sigma_1, \sigma_3, \sigma_4$, and $\sigma_5$ for our example path $\pi$.
\begin{figure}[ht]
\begin{center}
\includegraphics[scale=.6]{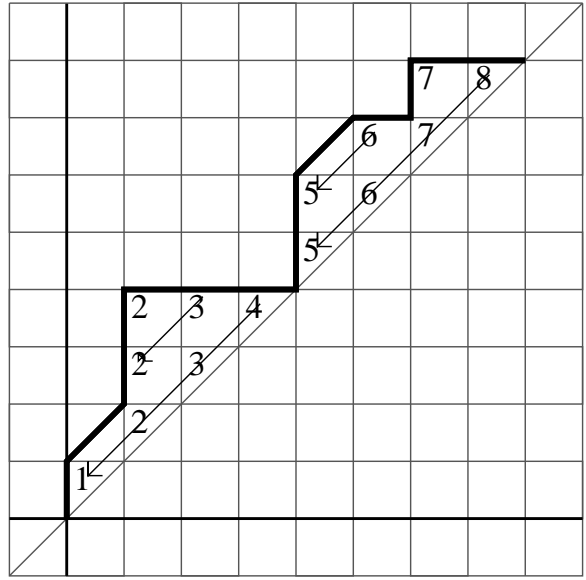}
\caption{The Schr\"oder path $\pi$ with $\sigma_1$, $\sigma_3$, 
$\sigma_4$, and $\sigma_5$ indicated.}
\end{center}
\end{figure}
In this example we have $\sigma_1 = s_8 s_7 s_6 s_5, \sigma_2 = s_7, \sigma_3 = s_6 s_5, \sigma_4 = s_4 s_3 s_2 s_1, \sigma_5 = s_3 s_2,$ and $\sigma_6 = s_2$.

We this set-up in mind, we can now define $\varphi(\pi)$.

\begin{definition}
\label{defn:varphi}
For any Schr\"oder path $\pi$, let $\sigma_1, \sigma_2, \ldots, \sigma_k$ denote the maps determined as in the discussion above.
Then we write
$$\varphi(\pi) = \sigma_k \sigma_{k-1} \ldots \sigma_1 (n+1,n,n-1, \ldots, 3,2,1).$$
\end{definition}

Summarizing our running example, we have the following.

\begin{example}
\label{exmp:running}
Let $\pi$ denote the Schr\"oder path in Figure 2 above, so that $\pi$ is given by $\pi = NDNNEEENNDENEE$.
Then $\sigma_1 = s_8 s_7 s_6 s_5, \sigma_2 = s_7, \sigma_3 = s_6 s_5, \sigma_4 = s_4 s_3 s_2 s_1, \sigma_5 = s_3 s_2, \sigma_6 = s_2$, and $\varphi(\pi) = 836791425$.
\end{example}

\begin{example}
Let $\pi$ denote the Schr\"oder path given by $\pi = DNEDDNNENDEE$.
Then $\sigma_1 = s_8 s_7 s_6 s_5$, $\sigma_2 = s_7 s_6$, $\sigma_3 = s_5$, $\sigma_4 = s_2$, and $\varphi(\pi) = 978624135$.
\end{example}

\begin{example}
Let $\pi$ denote the Schr\"oder path given by $\pi = NNENDNNEDEEDE$.
Then $\sigma_1 = s_8 s_7 s_6 s_5 s_4 s_3 s_2 s_1$, $\sigma_2 = s_6 s_5 s_4 s_3 s_2$, $\sigma_3 = s_5 s_4 s_3$, $\sigma_4 = s_3$, $\sigma_5 = s_1$, and $\varphi(\pi) = 683425719$.
\end{example}

Our goal for the remainder of this section is to show that the map $\varphi$ is a bijection which preserves the statistics $\tau_k$ for all $k \ge 1$.
To do this, we first consider the value of $\varphi$ on a Schr\"oder path which begins with a diagonal step.

\begin{proposition}
For all $n \ge 0$ and all $\pi \in \cS_n$ we have
\label{prop:varphiD}
\begin{equation}
\label{eqn:varphiD}
\varphi(D, \pi) = n+1, \varphi(\pi).
\end{equation}
\end{proposition}
\begin{proof}
Observe that by the construction of $\varphi(D,\pi)$ no $s_1$ will appear in $\sigma_k \sigma_{k-1} \ldots \sigma_1$.
Therefore the left-most element of $n+1, n, n-1, \ldots, 2,1$ will not be moved by $\sigma_k \sigma_{k-1} \ldots \sigma_1$.
Now the result follows.
\end{proof}

Next we consider the value of $\varphi$ on a product of two Schr\"oder paths.

\begin{proposition}
For all Schr\"oder paths $\pi_1$ and $\pi_2$ we have
\begin{equation}
\label{eqn:varphistar}
\varphi(\pi_1 * \pi_2) = \varphi(\pi_1) * \varphi(\pi_2).
\end{equation}
\end{proposition}
\begin{proof}
In the construction of $\varphi(\pi_1 * \pi_n)$, let $\sigma_1, \ldots, \sigma_j$ denote the strings of $s_i$s constructed from triangles below $\pi_2$, and observe that the next string constructed will be $s_i s_{i-1} \ldots s_1$, where $\pi_1$ ends at $(i-1, i-1)$.
Now let $\sigma_{j+2}, \ldots \sigma_k$ denote the remaining strings and let $a$ denote the left-most element of $\varphi(\pi_2)$.
Then we have
\begin{eqnarray*}
\sigma(\pi_1 * \pi_2) &=& \sigma_k \ldots \sigma_{j+2} s_i s_{i-1} \ldots s_1 \sigma_j \ldots \sigma_1 (n+1, n, \ldots, 2,1) \\[2ex]
&=& \sigma_k \ldots \sigma_{j+2} s_i s_{i-1} \ldots s_1 (n+1,n,\ldots, n-i+1, a, \tilde{\varphi(\pi_2)}) \\[2ex]
&=& \sigma_k \ldots \sigma_{j+2} (n, n-1, \ldots, n-i+1, a, n+1, \tilde{\varphi(\pi_2)}) \\[2ex]
&=& \tilde{\varphi(\pi_1)}, n+1, \tilde{\varphi(\pi_2)} \\[2ex]
&=& \varphi(\pi_1) * \varphi(\pi_2),
\end{eqnarray*}
as desired.
\end{proof}

We now show that $\varphi$ is a bijection from $\cS_n$ to $\fS_{n+1}(1243, 2143)$.

\begin{proposition}
\label{prop:varphibijection}
\renewcommand\labelenumi{{\upshape (\roman{enumi}) }}
\begin{enumerate}
\item
For all $\pi \in \cS$, we have $\varphi(\pi) \in \fS_{n+1}(1243, 2143)$.
\item
For all $n \ge 0$ the map $\varphi : \cS_n \longrightarrow \fS_{n+1}(1243, 2143)$ is a bijection.
\end{enumerate}
\end{proposition}
\begin{proof}
(i)
Observe that $\varphi(\emptyset) = 1$, $\varphi(D) = 21$, and $\varphi(NE) = 12$, so (i) holds for all $\pi \in \cS_0$ and all $\pi \in \cS_1$.
Arguing by induction, suppose (i) holds for all $\pi \in \cS_i$, where $0 \le i < n$, and fix $\pi \in \cS_n$.
If $\pi$ begins with a diagonal step then $\pi = D,\pi_1$ and by (\ref{eqn:varphiD}) we have $\varphi(\pi) = n+1, \varphi(\pi_1) \in \fS_{n+1}(1243, 2143)$.
If $\pi$ does not begin with a diagonal step then by Proposition \ref{prop:pathinirecurrence}(i) there exist $\pi_1 \in \cS_{i-1}$ and $\pi_2 \in \cS_{n-i}$, where $1 \le i \le n$, such that $\pi = \pi_1 * \pi_2$.
Then by (\ref{eqn:varphistar}) we have $\varphi(\pi) = \varphi(\pi_1) * \varphi(\pi_2) \in \fS_{n+1}(1243, 2143)$.
Now (i) follows.

(ii)
First observe that by (i) and since $|\cS_n| = |\fS_{n+1}(1243, 2143)|$ for all $n \ge 0$, it is sufficient to show that $\varphi$ is surjective.
To do this, we argue by induction on $n$.

It is routine to verify that $\varphi$ is surjective for $n = 0$ and $n = 1$, so suppose by induction that $\varphi$ is surjective for $1, 2, \ldots, n-1$ and fix $\pi \in \fS_{n+1}(1243, 2143)$.
If $\pi(1) = n+1$ then by Proposition \ref{prop:permnrecurrence}(i) there exists $\pi_1 \in \fS_n(1243, 2143)$ such that $\pi = n+1, \pi_1$.
By induction there exists $\alpha_1 \in \cS_{n-1}$ such that $\varphi(\alpha_1) = \pi_1$.
We now have
$$
\begin{array}{rclll}
\varphi(D, \alpha_1) &=& n+1, \pi_1 & \hspace{30pt} & \hbox{(by (\ref{eqn:varphiD}))} \\[2ex]
&=& \pi. \\
\end{array}
$$
If $\pi(1) \neq n+1$ then there exists $i$, $1 \le i \le n-1$, such that $\pi(i+1) = n+1$.
By Proposition \ref{prop:perminirecurrence}(i) there exist $\pi_1 \in \fS_i(1243, 2143)$ and $\pi_2 \in \fS_{n-i+1}(1243, 2143)$ such that $\pi = \pi_1 * \pi_2$.
By induction there exist $\alpha_1 \in \cS_{i-1}$ and $\alpha_2 \in \cS_{n-i}$ such that $\varphi(\alpha_1) = \pi_1$ and $\varphi(\alpha_2) = \pi_2$.
We now have
$$
\begin{array}{rclll}
\varphi(\alpha_1 * \alpha_2) &=& \pi_1 * \pi_2 &\hspace{30pt}& \hbox{(by (\ref{eqn:varphistar}))} \\[2ex]
&=& \pi.
\end{array}
$$
It follows that $\varphi$ is surjective, as desired.
\end{proof}

We now show that $\varphi$ preserves the statistics $\tau_k$ for all $k \ge 1$.

\begin{proposition}
\label{prop:varphitau}
For any $k \ge 1$ and any Schr\"oder path $\pi$ we have
\begin{equation}
\label{eqn:varphitau}
\tau_k(\varphi(\pi)) = \tau_k(\pi).
\end{equation}
\end{proposition}
\begin{proof}
Suppose $\pi \in \cS_n$;  we argue by induction on $n$.
The cases $n = 0$ and $n = 1$ are routine to verify, so suppose the result holds for all $\pi \in \cS_i$, where $i \le n-1$, and fix $\pi \in \cS_n$.
If $\pi$ begins with a diagonal step then $\pi = D, \pi_1$ for some $\pi_1 \in \cS_{n-1}$ and we have
$$
\begin{array}{rclll}
\tau_k(\varphi(\pi)) &=& \tau_k(\varphi(D, \pi_1)) & \hspace{30pt} & \\[2ex]
&=& \tau_k(n+1, \varphi(\pi_1)) & & \hbox{(by (\ref{eqn:varphiD}))} \\[2ex]
&=& \delta_{k1} + \tau_k(\varphi(\pi_1)) & & \hbox{(by Proposition \ref{prop:permnrecurrence}(ii),(iii))} \\[2ex]
&=& \delta_{k1} + \tau_k(\pi_1) & & \hbox{(by induction)} \\[2ex]
&=& \tau_k(D,\pi_1) & & \hbox{(by Proposition \ref{prop:pathnrecurrence}(ii),(iii))} \\[2ex]
&=& \tau_k(\pi).
\end{array}
$$
If $\pi$ does not begin with a diagonal step then $\pi = \pi_1 * \pi_2$ for some $\pi_1 \in \cS_{i-1}$ and $\pi_2 \in \cS_{n-i}$ and we have
$$
\begin{array}{rclll}
\tau_k(\varphi(\pi)) &=& \tau_k(\varphi(\pi_1 * \pi_2)) & \hspace{30pt} & \hbox{(by Proposition \ref{prop:pathinirecurrence}(i))} \\[2ex]
&=& \tau_k(\varphi(\pi_1) * \varphi(\pi_2)) & & \hbox{(by (\ref{eqn:varphistar}))} \\[2ex]
&=& \tau_k(\varphi(\pi_1)) + \tau_{k-1}(\varphi(\pi_1)) + \tau_k(\varphi(\pi_2)) & & \hbox{(by (\ref{eqn:tauperm}))} \\[2ex]
&=& \tau_k(\pi_1) + \tau_{k-1}(\pi_1) + \tau_k(\pi_2) & & \hbox{(by induction)} \\[2ex]
&=& \tau_k(\pi_1 * \pi_2) & & \hbox{(by (\ref{eqn:taupath}))} \\[2ex]
\end{array}
$$
as desired.
\end{proof}

In view of Proposition \ref{prop:varphitau}, the map $\varphi$ also relates the Schr\"oder paths of height at most $k-2$ with Schr\"oder permutations which avoid $12\ldots k$.
In particular, we have the following result.

\begin{corollary}
\label{cor:varphiheight}
Fix $k \ge 2$ and let $\cS_{n,k}$ denote the set of Schr\"oder paths in $\cS_n$ which do not cross the line $y - x = k-2$.
Then the restriction of $\varphi$ to $\cS_{n,k}$ is a bijection between $\cS_{n,k}$ and $\fS_{n+1}(1243, 2143, 12\ldots k)$.
\end{corollary}
\begin{proof}
Observe that $\pi \in \fS_{n+1}(1243, 2143)$ avoids $12\ldots k$ if and only if $\tau_k(\pi) = 0$.
By (\ref{eqn:varphitau}) this occurs if and only if $\tau_k(\varphi^{-1}(\pi)) = 0$, which occurs if and only if $\varphi^{-1}(\pi)$ does not cross the line $y - x = k-2$.
\end{proof}

For all $n \ge 0$, let $\cD_n$ denote the set of all Schr\"oder paths from $(0,0)$ to $(n,n)$ which contain no diagonal steps.
(Such paths are sometimes called {\it Catalan paths}, because $|\cD_n|$ is the well-known Catalan number $C_n = \frac{1}{n+1} {{2n} \choose {n}}$ for all $n \ge 0$.)
In \cite{Krattenthaler} Krattenthaler gives a bijection $\phi : \fS_n(132) \longrightarrow  \cD_n$ such that 
$$\tau_k(\pi) = \tau_k(\phi(\pi))$$
for all $k \ge 1$, all $n \ge 0$, and all $\pi \in \fS_n(132)$.
Krattenthaler's bijection $\phi$ is closely related to our bijection $\varphi$.
For any permutation $\pi$, let $\hat{\pi}$ denote the sequence obtained by adding one to every entry in $\pi$.
Observe that the map $\omega : \fS_n(132) \longrightarrow \fS_{n+1}(1243, 2143)$ given by $\omega(\pi) = 1, \hat{\pi}$ for all $\pi \in \fS_n(132)$ is a bijection between $\fS_n(132)$ and the set of permutations in $\fS_{n+1}(1243, 2143)$ which begin with 1.
The bijections $\phi$ and $\varphi$ are related in that
$$\phi(\pi) = \varphi^{-1}(\omega(\pi))$$
for all $n \ge 0$ and all $\pi \in \fS_n(132)$.

\section{The Continued Fractions}

In this section we will encounter several continued fractions, for which we will use the following notation.

\begin{definition}
For any given expressions $a_i$ $(i \ge 0)$ and $b_i$ $(i \ge 0)$ we write
$$\frac{a_0}{b_0} \cplus \frac{a_1}{b_1} \cplus \frac{a_2}{b_2} \cplus \frac{a_3}{b_3} \cplus \ldots$$
to denote the infinite continued fraction
$$
\frac{a_0}{\displaystyle
b_0 + \frac{a_1}{\displaystyle
b_1 + \frac{a_2}{\displaystyle
b_2 + \frac{a_3}{\displaystyle 
b_3 + \frac{a_4}{\displaystyle
b_4 + \cdots}}}}}.
$$
We use the corresponding notation for finite continued fractions.
\end{definition}

Several authors \cite{BCS,JR,Krattenthaler,MansourVainshtein4,MansourVainshtein5,RWZ} have described how to express the generating function for $\fS_n(132)$ with respect to various $\tau_k$ ($k \ge 1$) as a continued fraction.
The most general of these results is the following, which appears explicitly as \cite[Theorem 1]{BCS}.
A special case of this result appears as \cite[Theorem 1]{RWZ}, the result is implicit in \cite[Proposition 2.3]{MansourVainshtein4}, and it can be proved by modifying slightly the techniques of \cite[Corollary 7]{JR} and \cite[Theorem 1]{Krattenthaler}.

\begin{theorem}
\label{thm:132cf}
(Br\"and\'en, Claesson, and Steingr\'imsson \cite[Theorem 1]{BCS})
For all $i \ge 1$, let $x_i$ denote an indeterminate.
Then we have
$$\sum_{\pi \in \fS(132)} \prod_{k \ge 1} x_k^{\tau_k(\pi)} = \frac{1}{1} \cminus \frac{x_1}{1} \cminus \frac{x_1 x_2}{1} \cminus \frac{x_1 x_2^2 x_3}{1} \cminus \ldots \cminus \frac{\displaystyle \prod_{k \ge 1} x_k^{{{n} \choose {k}}}}{1} \cminus \ldots.$$
\end{theorem}

In the same paper, Br\"and\'en, Claesson, and Steingr\'imsson define a {\it Catalan continued fraction} to be a continued fraction of the form
$$\frac{1}{1} \cminus \frac{m_0}{1} \cminus \frac{m_1}{1} \cminus \frac{m_2}{1} \cminus \ldots,$$
where for all $i \ge 0$, the expression $m_i$ is a monic monomial in a given set of variables.
Roughly speaking, Br\"and\'en, Claesson, and Steingr\'imsson show \cite[Theorem 2]{BCS} that the multivariate generating function for $\fS(132)$ with respect to a given countable family of statistics may be expressed as a Catalan continued fraction if and only if each statistic in the family is a (possibly infinite) linear combination of the $\tau_k$s and each $\tau_k$ appears in only finitely many of these linear combinations.

In this section we prove analogues of these results for permutations which avoid 1243 and 2143.
We begin by adapting Krattenthaler's approach in \cite{Krattenthaler}, using our bijection $\varphi$ and a result of Flajolet to express the generating function for $\fS(1243, 2143)$ with respect to $\tau_k$, $k \ge 1$, as a continued fraction.
For convenience, we first recall the relevant specialization of Flajolet's result.

\begin{theorem}
\label{thm:Flajolet}
(Flajolet \cite[Theorem 1]{Flajolet})
For all $i \ge 1$, let $x_i$ denote an indeterminate.
Then we have
\begin{equation}
\label{eqn:Flajolet}
\sum_{\pi \in \cS} \prod_{k \ge 1} x_k^{\tau_k(\pi)} = 
\frac{x_1}{1 - x_1} \cminus \frac{x_1 x_2}{1 - x_1 x_2} \cminus \frac{x_1 x_2^2 x_3}{1 - x_1 x_2^2 x_3} \cminus \ldots \cminus \frac{\displaystyle \prod_{i=0}^n x_i^{{{n} \choose {i}}}}{\displaystyle 1 - \prod_{i=0}^n x_i^{{{n} \choose {i}}}} \cminus \ldots.
\end{equation}
\end{theorem}

Combining $\varphi$ with Theorem \ref{thm:Flajolet}, we obtain the following analogue of Theorem \ref{thm:132cf}.

\begin{theorem}
\label{thm:contfrac}
For all $i \ge 1$, let $x_i$ denote an indeterminate.
Then we have
\begin{equation}
\label{eqn:contfrac}
\sum_{\pi \in \fS(1243, 2143)} \prod_{k \ge 1} x_k^{\tau_k(\pi)} = 
1 + \frac{x_1}{1 - x_1} \cminus \frac{x_1 x_2}{1 - x_1 x_2} \cminus \frac{x_1 x_2^2 x_3}{1 - x_1 x_2^2 x_3} \cminus \ldots \cminus \frac{\displaystyle \prod_{i=0}^n x_i^{{{n} \choose {i}}}}{\displaystyle 1 - \prod_{i=0}^n x_i^{{{n} \choose {i}}}} \cminus \ldots.
\end{equation}
\end{theorem}
\begin{proof}
This is immediate from Theorem \ref{thm:Flajolet}, in view of Propositions \ref{prop:varphibijection} and \ref{prop:varphitau}.
\end{proof}

Using (\ref{eqn:contfrac}), we can express the generating function for $|\fS_n(1243, 2143, 12\ldots k)|$ as a (finite) continued fraction.

\begin{corollary}
\label{cor:12k}
For all $k \ge 1$ we have
\begin{equation}
\label{eqn:12gf}
\sum_{n=0}^\infty |\fS_n(1243, 2143, 12\ldots k)| x^n = 
1 + \underbrace{\frac{x}{1-x} \cminus \frac{x}{1-x} \cminus \ldots \cminus \frac{x}{1-x}}_{k-1\ terms}.
\end{equation}
\end{corollary}
\begin{proof}
In (\ref{eqn:contfrac}) set $x_1 = x$, $x_2 = x_3 = \ldots = x_{k-1} = 1$ and $x_i = 0$ for all $i \ge k$.
\end{proof}

Curiously, as we prove in Proposition \ref{prop:sporadfrac}, the continued fraction on the right side of (\ref{eqn:12gf}) is also equal to the generating function for $|\fS_n(1243, 2143, 213\ldots k)|$.

Using (\ref{eqn:contfrac}), we can also express the generating function for $\fS(1243, 2143)$ with respect to the total number of increasing subsequences as a continued fraction.

\begin{corollary}
For any permutation $\pi$, let $m(\pi)$ denote the number of nonempty increasing subsequences in $\pi$.
Then 
$$\sum_{\pi \in \fS(1243, 2143)} q^{m(\pi)} x^{|\pi|} = 1 + \frac{x q}{1 - xq} \cminus \frac{x q^2}{1 - xq^2} \cminus \frac{x q^4}{1 - xq^4} \cminus \ldots \cminus \frac{x q^{2^n}}{1 - x q^{2^n}} \cminus \ldots.$$ 
\end{corollary}
\begin{proof}
In (\ref{eqn:contfrac}), set $x_1 = xq$ and $x_i = q$ for all $i \ge 2$ and use the fact that $m = \sum\limits_{k \ge 1} \tau_k.$
\end{proof}

By specializing the $x_i$s in (\ref{eqn:contfrac}) in various ways, one can obtain continued fractions of many other interesting forms.
As an example, we have the following result.

\begin{corollary}
For any permutation $\pi$, let $m(\pi)$ denote the length of $\pi$ plus the number of noninversions in $\pi$.
Then
$$\sum_{\pi \in \fS(1243, 2143)} q^{m(\pi)} = 1 + \frac{q}{1-q} \cminus \frac{q^2}{1-q^2} \cminus \frac{q^3}{1 - q^3} \cminus \frac{q^4}{1 - q^4} \cminus \ldots.$$
\end{corollary}
\begin{proof}
In (\ref{eqn:contfrac}), set $x_1 = x_2 = q$ and $x_i = 1$ for all $i \ge 3$ and use the fact that $m = \tau_1 + \tau_2$.
\end{proof}

We now turn our attention to the question of which statistics on $\fS(1243, 2143)$ have generating functions which can be expressed as continued fractions like the one in (\ref{eqn:contfrac}). 
We begin by specifying which continued fractions we wish to consider.

By a {\it Schr\"oder continued fraction} we mean a continued fraction of the form
$$1 + \frac{m_0}{1 - m_0} \cminus \frac{m_1}{1 - m_1} \cminus \frac{m_2}{1 - m_2} \cminus \ldots,$$
in which $m_i$ is a monic monomial in a given set of variables for all $i \ge 0$.
Observe that if $f_1, f_2, f_3, \ldots$ are (possibly infinite) linear combinations of the $\tau_k$s with the property that each $\tau_k$ appears in only finitely many $f_i$, then by specializing the $x_i$s appropriately in (\ref{eqn:contfrac}) we can express the generating function
$$\sum_{\pi \in \fS(1243, 2143)} x^{|\pi|} \prod_{k \ge 1} q_k^{f_k(\pi)}$$
as a Schr\"oder continued fraction.
For instance, when only $f_1$ is present, we have the following corollary of Theorem \ref{thm:E1}.

\begin{corollary}
Let $\lambda_1, \lambda_2, \ldots$ denote nonnegative integers and let $f$ denote the statistic
$$f = \sum_{k \ge 1} \lambda_k \tau_k$$
on $\fS(1243, 2143)$.
Then
\begin{eqnarray*}
\lefteqn{\sum_{\pi \in \fS(1243, 2143)} q^{f(\pi)} x^{|\pi|} =} & & \\
& & 1 + \frac{x q^{f(1)}}{1 - x q^{f(1)}} \cminus \frac{x q^{f(12) - f(1)}}{1 - x q^{f(12) - f(1)}} \cminus \frac{x q^{f(123) - f(12)}}{1 - x q^{f(123) - f(12)}} \cminus \frac{x q^{f(1234) - f(123)}}{1 - x q^{f(1234) - f(123)}} \cminus \ldots. \\
\end{eqnarray*}
\end{corollary}
\begin{proof}
In (\ref{eqn:contfrac}) set $x_1 = x q^{\lambda_1}$ and $x_i = q^{\lambda_i}$ for all $i \ge 2$ to obtain
$$\sum_{\pi \in \fS(1243, 2143)} q^{f(\pi)} x^{|\pi|} = 1 + \frac{xq^{\lambda_1}}{1- x q^{\lambda_1}} \cminus \frac{x q^{\lambda_1 + \lambda_2}}{1 - x q^{\lambda_1 + \lambda_2}} \cminus \frac{x q^{\lambda_1 + 2 \lambda_2 + \lambda_3}}{1 - x q^{\lambda_1 + 2 \lambda_2 + \lambda_3}} \cminus \ldots.$$
Now the result follows from the fact that $f(123\ldots k) - f(12 \ldots k-1) = \sum\limits_{i=0}^{k-1} {{k-1} \choose {i}} \lambda_i$ for all $k \ge 2$.
\end{proof}

Modifying the proof of \cite[Theorem 2]{BCS} slightly, we now show that linear combinations of the $\tau_k$s are the only statistics on $\fS(1243, 2143)$ whose generating functions can be expressed as Schr\"oder continued fractions.
In order to do this, we first set some notation.

Let $\cA$ denote the ring of matrices with integer entries whose rows and columns are indexed by the positive integers, and which have the property that in any row only finitely many entries are nonzero.
Then each element of $\cA$ corresponds to an infinite family of statistics on $\fS(1243, 2143)$.

\begin{definition}
\label{defn:Agf}
For all $A \in \cA$ and all $n \ge 1$, let $\tau_{A,n}$ denote the linear combination of the $\tau_k$s whose coefficients appear in the $n$th column of $A$.
That is,
\begin{equation}
\label{eqn:tauAn}
\tau_{A,n} = \sum_{k \ge 1} A_{kn} \tau_k
\end{equation}
for all $A \in \cA$ and all $n \ge 1$.
\end{definition}

Since each element of $\cA$ corresponds to a family of statistics on $\fS(1243, 2143)$, each element also has an associated multivariate generating function.

\begin{definition}
\label{defn:Acf}
For all $i \ge 1$, let $q_i$ denote an indeterminate.
For all $A \in \cA$, we write $F_A({\bf q})$ to denote the generating function given by
\begin{equation}
\label{eqn:Agf}
F_A({\bf q}) = \sum_{\pi \in \fS(1243, 2143)} \prod_{k \ge 1} q_k^{\tau_{A,k}(\pi)}.
\end{equation}
\end{definition}

We also associate with each element of $\cA$ a Schr\"oder continued fraction.

\begin{definition}
For all $i \ge 1$, let $q_i$ denote an indeterminate.
For all $A \in \cA$, we write $C_A({\bf q})$ to denote the continued fraction given by
\begin{equation}
\label{eqn:Acf}
C_A({\bf q}) = 1 + \frac{\prod\limits_{k \ge 1} q_k^{A_{1k}}}{1 - \prod\limits_{k \ge 1} q_k^{A_{1k}}} \cminus \frac{\prod\limits_{k \ge 1} q_k^{A_{2k}}}{1 - \prod\limits_{k \ge 1} q_k^{A_{2k}}} \cminus \frac{\prod\limits_{k \ge 1} q_k^{A_{3k}}}{1 - \prod\limits_{k \ge 1} q_k^{A_{3k}}} \cminus \ldots.
\end{equation}
\end{definition}

We now give our analogue of \cite[Theorem 2]{BCS}.
Although the proof is nearly identical to the proof of \cite[Theorem 2]{BCS}, we include it here for completeness.

\begin{theorem}
\label{thm:BCS}
Let $B$ denote the matrix in $\cA$ which satisfies $B_{ij} = {{i-1} \choose {j-1}}$.
(Here we use the convention that ${{i} \choose {j}} = 0$ whenever $i < j$.)
Then for all $A \in \cA$,
\begin{equation}
\label{eqn:FACBA}
F_A({\bf q}) = C_{BA}({\bf q})
\end{equation}
and
\begin{equation}
\label{eqn:CAFB-1A}
C_A({\bf q}) = F_{B^{-1}A}({\bf q}).
\end{equation}
In particular, the set of Schr\"oder continued fractions is exactly the set of generating functions for countable families of statistics on $\fS(1243, 2143)$ in which each statistic is a (possibly infinite) linear combination of the $\tau_k$s and each $\tau_k$ appears in only finitely many statistics.
\end{theorem}
\begin{proof}
To prove (\ref{eqn:FACBA}), we apply (\ref{eqn:contfrac}) with $x_i = \prod\limits_{j \ge 1} q_j^{A_{ij}}$ for all $i \ge 1$.
We have
$$
\begin{array}{rclll}
C_{BA}({\bf q}) &=& {\displaystyle 1 + \frac{x_1}{1 - x_1} \cminus \frac{x_1 x_2}{1 - x_1 x_2} \cminus \frac{x_1 x_2^2 x_3}{1 - x_1 x_2^2 x_3} \cminus \ldots} & \hspace{-1pt} & \hbox{(by (\ref{eqn:Acf}))} \\[4ex]
&=& {\displaystyle \sum_{\pi \in \fS(1243, 2143)} \prod\limits_{k \ge 1} x_k^{\tau_k(\pi)}} & & \hbox{(by (\ref{eqn:contfrac}))} \\[2ex]
&=& {\displaystyle \sum_{\pi \in \fS(1243, 2143)} \prod\limits_{k \ge 1} \left( \prod\limits_{j \ge 1} q_j^{A_{kj}}\right)^{\tau_k(\pi)}} & & \\[4ex]
&=& {\displaystyle \sum_{\pi \in \fS(1243, 2143)} \prod\limits_{j \ge 1} q_j^{\sum\limits_{k\ge 1} A_{kj} \tau_k(\pi)}} & & \\[4ex]
&=& {\displaystyle F_A({\bf q})} & & \hbox{(by (\ref{eqn:tauAn}) and (\ref{eqn:Agf})),}   
\end{array}
$$
as desired.
To prove (\ref{eqn:CAFB-1A}), observe that $(B^{-1})_{ij} = (-1)^{i+j} {{i-1} \choose {j-1}}$ so $B^{-1} \in \cA$.
Therefore we may replace $A$ with $B^{-1} A$ in (\ref{eqn:FACBA}) to obtain (\ref{eqn:CAFB-1A}). 
\end{proof}

\noindent
{\bf Remark}
There is a result for finite Schr\"oder continued fractions which is analogous to Theorem \ref{thm:BCS}.
Specifically, fix $k \ge 2$ and let $\cA_k$ denote the ring of $k-1$ by $k-1$ matrices with integer entries.
If we replace $\cA$ with $\cA_k$ and $\fS(1243, 2143)$ with $\fS(1243, 2143, 12\ldots k)$ in Definitions \ref{defn:Agf} and \ref{defn:Agf} then the resulting analogue of Theorem \ref{thm:BCS} follows by an argument almost identical to the proof of Theorem \ref{thm:BCS}.

\section{Generating Functions Involving Chebyshev Polynomials}
\label{sec:Chebyshev}

In this section we will be interested in connections between pattern-avoiding permutations and Chebyshev polynomials of the second kind.
We begin by recalling these polynomials.

\begin{definition}
For all $n \ge -1$, we write $U_n(x)$ to denote the {\em $n$th Chebyshev polynomial of the second kind}, which is defined by $U_{-1}(x) = 0$ and ${\displaystyle U_n(\cos t) = \frac{\sin((n+1) t)}{\sin t}}$ for $n \ge 0$.
These polynomials satisfy
\begin{equation}
\label{eqn:Unrecurrence}
U_n(x) = 2 x U_{n-1}(x) - U_{n-2}(x) \hspace{30pt} (n \ge 1).
\end{equation}
\end{definition}

We will find it useful to reformulate the recurrence in (\ref{eqn:Unrecurrence}) as follows.

\begin{lemma}
For all $n \ge 1$,
\begin{equation}
\label{eqn:qnUnrecurrence}
U_n\left(\frac{1-x}{2\sqrt{x}}\right) = \frac{1-x}{\sqrt{x}} U_{n-1}\left(\frac{1-x}{2\sqrt{x}}\right) - U_{n-2}\left(\frac{1-x}{2\sqrt{x}}\right).
\end{equation}
\end{lemma}
\begin{proof}
This is immediate from (\ref{eqn:Unrecurrence}).
\end{proof}

Chebyshev polynomials of the second kind have been found to have close connections with $\fS(132)$;  see, for instance, \cite{CW,Krattenthaler,MansourVainshtein4,MansourVainshtein,MansourVainshtein5}.
Two such connections are given in the following results.

\begin{theorem}
\label{thm:C1}
(Krattenthaler \cite[Theorem 2]{Krattenthaler}, Mansour and Vainshtein \cite[Theorem 3.1]{MansourVainshtein4}, Chow and West \cite[Theorem 3.6, second case]{CW})
For all $k \ge 1$ we have
$$\sum_{n=0}^\infty |\fS_n(132, 12\ldots k)| x^n = \frac{U_{k-1}\left( \frac{1}{2 \sqrt{x}} \right)}{\sqrt{x} U_k\left( \frac{1}{2 \sqrt{x}}\right)}.$$
\end{theorem}

\begin{theorem}
\label{thm:C2}
(Krattenthaler \cite[Theorem 3]{Krattenthaler}, Mansour and Vainshtein \cite[Theorems 3.1 and 4.1]{MansourVainshtein4})
Let $r \ge 1$, $b \ge 1$, and $k \ge 2$ satisfy ${{k+b-1}\choose{k}} \le r < {{k+b} \choose {k}}$.
Then the generating function for the number of 132-avoiding permutations which contain exactly $r$ subsequences of type $12\ldots k$ is given by 
$$\sum \prod_{i = 2}^b {{l_{i-1} + l_i - 1}\choose {l_i}} \frac{\left( U_{k-1}\left(\frac{1}{2\sqrt{x}}\right)\right)^{l_1 - 1}}{\left(U_k\left( \frac{1}{2\sqrt{x}}\right)\right)^{l_1 + 1}} x^{\frac{1}{2}( l_1 - 1 ) + \sum\limits_{j=2}^b l_j},$$
where the sum on the left is over all sequences $l_1, l_2, \ldots, l_b$ of nonnegative integers such that
$$\sum_{i=1}^b l_i {{k+i-2}\choose{k-1}} = r.$$
\end{theorem}

In this section we prove analogues of Theorems \ref{thm:C1} and \ref{thm:C2} for $\fS_n(1243, 2143)$.
We begin with an analogue of Theorem \ref{thm:C1}.

\begin{theorem}
\label{thm:E1}
For all $k \ge 1$ we have
\begin{equation}
\label{eqn:E112k}
\sum_{n=0}^\infty |\fS_n(1243, 2143, 12\ldots k)| x^n = 1 + \frac{\sqrt{x} U_{k-2}\left(\frac{1-x}{2 \sqrt{x}}\right)}{U_{k-1}\left( \frac{1-x}{2 \sqrt{x}}\right)}
\end{equation}
and
\begin{equation}
\label{eqn:E121k}
\sum_{n=0}^\infty |\fS_n(1243, 2143, 213\ldots k)| x^n = 1 + \frac{\sqrt{x} U_{k-2}\left(\frac{1-x}{2 \sqrt{x}}\right)}{U_{k-1}\left( \frac{1-x}{2 \sqrt{x}}\right)}.
\end{equation}
\end{theorem}

Before pressing on to the proof of Theorem \ref{thm:E1}, we observe that when we set $k = 3$ in (\ref{eqn:E112k}) and (\ref{eqn:E121k}) we find that
$$\sum_{n=0}^\infty |\fS_n(2143, 123)|x^n = 1 + x \frac{1-x}{1-3x+x^2}$$
and
$$\sum_{n=0}^\infty |\fS_n(1243, 213)|x^n = 1 + x \frac{1-x}{1-3x+x^2}.$$
These generating functions were originally found by West \cite[Example 9 and Table 1]{WestGenTrees}, using generating trees.

Our next result is the key to our proof of Theorem \ref{thm:E1}.
To state it, we first set some notation.

\begin{definition}
For any set $R$ of permutations, let $Rn$ denote the set of permutations obtained by appending the smallest possible positive integer to each permutation in $R$.
\end{definition}

\begin{example}
If $R = \{42531, 4231, 312\}$ then $Rn = \{425316, 42315, 3124\}$.
\end{example}

\begin{theorem}
\label{thm:pik}
Let $R$ denote a nonempty set of permutations and set
\begin{equation}
\label{eqn:P}
P(x) = \sum_{n=0}^\infty |\fS_n(1243, 2143, R)| x^n
\end{equation}
and
\begin{equation}
\label{eqn:Q}
Q(x) = \sum_{n=0}^\infty |\fS_n(1243, 2143, Rn)| x^n.
\end{equation}
Then
\begin{equation}
\label{eqn:QP}
Q(x) = \frac{2 - P(x)}{2 - x - P(x)}.
\end{equation}
\end{theorem}
\begin{proof}
Arguing as in the proof of Proposition \ref{prop:permnrecurrence}(i), we find that for all $n \ge 0$ the map
$$
\begin{array}{ccc}
\fS_n(1243, 2143, Rn) &\longrightarrow& \fS_{n+1}(1243, 2143, Rn) \\
\pi &\mapsto& n+1,\pi
\end{array}
$$
is a bijection between $\fS_n(1243, 2143, Rn)$ and the set of all permutations in $\fS_{n+1}(1243, 2143, Rn)$ which begin with $n+1$.
Arguing as in the proof of Proposition \ref{prop:perminirecurrence}(i), we find that for all $i$ and all $n$ such that $1 \le i \le n-1$, the map
$$
\begin{array}{ccc}
\fS_i(1243, 2143, R) \times \fS_{n-i}(1243, 2143, Rn) &\longrightarrow& \fS_n(1243, 2143, Rn) \\
(\pi_1, \pi_2) &\mapsto& \pi_1 * \pi_2
\end{array}
$$
is a bijection between $\fS_i(1243, 2143, R) \times \fS_{n-i}(1243, 2143, Rn)$ and the set of permutations in $\fS_n(1243, 2143, Rn)$ for which $\pi(i+1) = n$.
Combining these two bijections, we find that
$$Q(x) = 1 + x Q(x) + (P(x) - 1)(Q(x) - 1).$$
Solving this last equation for $Q(x)$, we obtain (\ref{eqn:QP}).
\end{proof}

As we show next, if $P(x)$ is a rational function then $Q(x)$ has a simple form.

\begin{corollary}
\label{cor:ratpik}
Let $R$ denote a nonempty set of permutations.
If
$$\sum_{n=0}^\infty |\fS_n(1243, 2143, R)|x^n = 1 + x \frac{f(x)}{g(x)}$$
then
$$\sum_{n=0}^\infty |\fS_n(1243, 2143, Rn)|x^n = 1 + x \frac{g(x)}{(1-x)g(x) - x f(x)}.$$
\end{corollary}
\begin{proof}
By (\ref{eqn:QP}) we have
\begin{eqnarray*}
\sum_{n=0}^\infty |\fS_n(1243, 2143, Rn)|x^n &=& \frac{\displaystyle 1 - x \frac{f(x)}{g(x)}}{\displaystyle 1 - x + x \frac{f(x)}{g(x)}} \\[3ex]
&=& 1 + x \frac{g(x)}{(1-x)g(x) - x f(x)},
\end{eqnarray*}
as desired.
\end{proof}

Using Corollary \ref{cor:ratpik} it is now routine to prove Theorem \ref{thm:E1}.

\bigskip
\noindent
{\it Proof of Theorem \ref{thm:E1}. }
We argue by induction on $k$.
It is routine to verify that (\ref{eqn:E112k}) holds for $k = 1$ and $k = 2$, so we assume (\ref{eqn:E112k}) holds for $k-2$ and $k-1$.
Then by Corollary \ref{cor:ratpik} and (\ref{eqn:qnUnrecurrence}) we have
\begin{eqnarray*}
\sum_{n=0}^\infty |\fS_n(1243, 2143, 12\ldots k)| x^n &=& 1 + \frac{\sqrt{x} U_{k-1}\left(\frac{1-x}{2\sqrt{x}}\right)}{\frac{1-x}{\sqrt{x}} U_{k-1}\left(\frac{1-x}{2\sqrt{x}}\right) - U_{k-2}\left(\frac{1-x}{2\sqrt{x}}\right)} \\[2ex]
&=& 1 + \frac{\sqrt{x} U_{k-2}\left(\frac{1-x}{2 \sqrt{x}}\right)}{U_{k-1}\left( \frac{1-x}{2 \sqrt{x}}\right)},
\end{eqnarray*}
as desired.
The proof of (\ref{eqn:E121k}) is similar to the proof of (\ref{eqn:E112k}).
$\Box$

\bigskip

We now have a corollary of Theorem \ref{thm:E1} involving the continued fraction which appears in (\ref{eqn:12gf}).

\begin{corollary}
\label{prop:sporadfrac}
For all $k \ge 2$ we have
$$\sum_{n=0}^\infty |\fS_n(1243, 2143, 2 1 3 4 \ldots k)| x^n = 
1 + \underbrace{\frac{x}{1-x} \cminus \frac{x}{1-x} \cminus \ldots \cminus \frac{x}{1-x}}_{k-1\ terms}.
$$
\end{corollary}
\begin{proof}
This is immediate from Theorem \ref{thm:E1} and Corollary \ref{cor:12k}.
\end{proof} 

Next we turn our attention to our analogue of Theorem \ref{thm:C2}.

\begin{theorem}
\label{thm:E2}
Fix $r \ge 1$, $k \ge 2$, and $b \ge 0$ such that ${{k+b} \choose {k}} \le r < {{k+b+1} \choose {k}}$.
Then we have
\begin{equation}
\label{eqn:12krgf}
\sum_{\pi} x^{|\pi|} = \sum \prod_{i=0}^b {{l_i + l_{i+1} + m_i - 1} \choose {l_{i+1} + m_i}} {{l_{i+1} + m_i} \choose {m_i}} \frac{\left(U_{k-2}\left(\frac{1-x}{2\sqrt{x}}\right)\right)^{l_0 - 1}}{\left( U_{k-1}\left(\frac{1-x}{2\sqrt{x}}\right)\right)^{l_0+1}} x^{\frac{1}{2}(1-l_0) + \sum\limits_{j=0}^b ( l_j + m_j )}.
\end{equation}
Here the sum on the left is over all permutations in $\fS(1243, 2143)$ which contain exactly $r$ subsequences of type $12 \ldots k$.
The sum on the right is over all sequences $l_0, l_1, \ldots, l_b,$ and $m_0, m_1, \ldots, m_b$ of nonnegative integers such that 
\begin{equation}
\label{eqn:rcondition}
r = \sum_{i=0}^b (l_i + m_i) {{k+i-1} \choose {k-1}}.
\end{equation}
Throughout we adopt the convention that ${{a} \choose {0}} = 1$ and ${{a} \choose {-1}} = 0$ for any integer $a$.
\end{theorem}

To prove this theorem, we first need to set some notation and prove some preliminary results.
We begin with a certain matrix.

\begin{definition}
For all $k \ge 0$, we write $A_k$ to denote the $k+1$ by $k+1$ tridiagonal matrix given by
$$A_k = \left(
\matrix{x & \sqrt{x} & 0 & 0 & 0 & \cdots & 0 & 0 & 0 & 0\cr
\sqrt{x} & x & \sqrt{x} & 0 & 0 & \cdots & 0 & 0 & 0 & 0 \cr
0 & \sqrt{x} & x & \sqrt{x} & 0 & \cdots & 0 & 0 & 0 & 0\cr
\vdots & & \ddots & \ddots & \ddots & & & & \vdots & \vdots \cr
\vdots & & & \ddots & \ddots & \ddots & & & \vdots & \vdots \cr
\vdots & & & & \ddots & \ddots & \ddots & & & \vdots \cr
\vdots & & & & & \ddots & \ddots & \ddots & \vdots & \vdots \cr
0 & 0 & 0 & 0 & 0 & \cdots & \sqrt{x} & x & \sqrt{x} & 0 \cr
0 & 0 & 0 & 0 & 0 & \cdots & 0 & \sqrt{x} & x & \sqrt{x} \cr
0 & 0 & 0 & 0 & 0 & \cdots & 0 & 0 & \sqrt{x} & x}
\right).$$
\end{definition}

The matrix $A_k$ is closely related to generating functions for various sets of Schr\"oder paths.
To describe this relationship, we let $\cS(r,s)$ denote the set of lattice paths involving only east $(1,0)$, north $(0,1)$, and diagonal $(1,1)$ steps which begin at a point at height $r$, end at a point at height $s$, and do not cross the lines $y - x = k$ and $y = x$.
For any such path $\pi$ which begins at $(x_1, y_1)$ and ends at $(x_2, y_2)$, we write $l(\pi) = \frac{1}{2}\left( x_2 + y_2 - x_1 - y_1 \right)$.
The next lemma summarizes the relationship between $A_k$ and the generating function for $\cS(r,s)$ with respect to $l$.

\begin{lemma}
\label{lem:pathheightgf}
For all $k \ge 0$ and all $r$ and $s$ such that $0 \le r,s \le k$ we have 
$$\sum_{\pi \in \cS(r,s)} x^{l(\pi)} = \frac{(-1)^{r+s} \det(I - A_k; s,r)}{\det(I - A_k)}.$$
Here $I$ is the identity matrix of the appropriate size and $\det(I - A_k; s,r)$ is the minor of $I - A_k$ in which the $s$th row and $r$th column have been deleted.
\end{lemma}

For the sake of brevity we omit the proof of Lemma \ref{lem:pathheightgf}.
An outline of this proof can be obtained by modifying slightly the outline of the proof of \cite[Theorem A2]{Krattenthaler}.

The matrix $A_k$ is also closely connected with Chebyshev polynomials of the second kind.
In particular, we have the following result.

\begin{lemma}
\label{lem:qkdet}
For all $k \ge 0$,
$$(\sqrt{x})^{k+1} U_{k+1}\left(\frac{1-x}{2\sqrt{x}}\right) = \det(I - A_k).$$
\end{lemma}
\begin{proof}
We argue by induction on $k$.
It is routine to verify that the result holds for $k = 0$ and $k = 1$, so we assume the result holds for $k-2$ and $k-1$.
By expanding $\det(I - A_k)$ along the bottom row, one can show that $(\sqrt{x})^{-(k+1)} \det(I - A_k)$ satisfies (\ref{eqn:qnUnrecurrence}).
By our induction assumption we find the result holds for $k$, as desired.
\end{proof}

As a warm-up for our proof of Theorem \ref{thm:E2}, we use Lemmas \ref{lem:pathheightgf} and \ref{lem:qkdet} to give another proof of (\ref{eqn:E112k}).

\bigskip

\noindent
{\it Another proof of (\ref{eqn:E112k}). }
First observe that if $F(x)$ is the generating function on the left side of (\ref{eqn:E112k}) and $G(x)$ is the generating function for Schr\"oder paths from $(0,0)$ to $(n,n)$ which do not exceed height $k-2$, then in view of (\ref{eqn:varphitau}) we have $F(x) = 1 + x G(x)$.
Combining Lemmas \ref{lem:pathheightgf} and \ref{lem:qkdet} we find that ${\displaystyle G(x) = \frac{U_{k-2}\left(\frac{1-x}{2 \sqrt{x}}\right)}{\sqrt{x} U_{k-1}\left( \frac{1-x}{2 \sqrt{x}}\right)}}$.
Now (\ref{eqn:E112k}) follows.
$\Box$

\bigskip

We now prove Theorem \ref{thm:E2}.

\bigskip

\noindent
{\it Proof of Theorem \ref{thm:E2}. }
First observe that if $F(x)$ is the generating function on the left side of (\ref{eqn:12krgf}) and $G(x)$ is the generating function for the set of Schr\"oder paths from $(0,0)$ to $(n,n)$ for which $\tau_k(\pi) = r$, then in view of (\ref{eqn:varphitau}) we have $F(x) = x G(x)$.
With this in mind, we find $G(x)$.

To compute $G(x)$, we observe that any Schr\"oder path $\pi$ with $\tau_k(\pi) = r$ can be constructed by the following procedure in exactly one way.
\begin{enumerate}
\item
Choose $l_0, l_1, \ldots, l_b,$ and $m_0, m_1, \ldots, m_b$ such that (\ref{eqn:rcondition}) holds.
Construct a sequence of east and diagonal steps which contains exactly $l_i$ east steps at height $k+i-1$ and $m_i$ diagonal steps at height $k+i-1$ for $0 \le i \le b$ and which satisfies both of the following.
\begin{enumerate}
\item The step immediately preceeding a step at height $j$ is either an east step at height $j+1$ or less, or a diagonal step at height $j$ or less.
\item All steps after the last east step at height $j$ are at height $j - 1$ or less.
\end{enumerate}

\item
After each east step at height $k-1$ except the last, insert an (possibly empty) upside-down Schr\"oder path of height at most $k-2$.

\item 
Before the first step insert a path from height 0 to height $k-2$ which does not exceed height $k-2$.

\item 
After the last step insert a path from height $k-2$ to height 0 which does not exceed height $k-2$.
\end{enumerate}
Since the choice at each step is independent of the choices at the other steps, and since every sequence of choices results in a path of the type desired, $G(x)$ is the product of the generating functions for each step.

To compute the generating function for step 1, suppose we have fixed $l_0, l_1, \ldots, l_b,$ and $m_0, m_1, \ldots, m_b$;  then each of the resulting partial paths will have generating function $x^{\sum\limits_{j=0}^b \left( l_i + m_i\right)}$.
To count these paths, we construct them from the top down.
That is, we first arrange the $m_b$ diagonal steps at height $b$;  there is one way to do this.
We then place the $l_b$ east steps at height $b$ so that one of these steps occurs after all of the diagonal steps.
There are ${{l_b + m_b - 1} \choose {m_b}}$ ways to do this.
We then place the $m_{b-1}$ diagonal steps at height $b-1$ so that none of these steps immediately follows a diagonal step at height $b$.
There are ${{l_b + m_{b-1}} \choose {m_{b-1}}}$ ways to do this.
Proceeding in this fashion, we find that the generating function for step 1 is
\begin{equation}
\label{eqn:step1}
\sum \prod_{i=0}^b {{l_i + l_{i+1} + m_i - 1} \choose {l_{i+1} + m_i}} {{l_{i+1} + m_i} \choose {m_i}} x^{\sum\limits_{j=0}^b \left( l_i + m_i\right)},
\end{equation}
where the sum is over all sequences $l_0, l_1, \ldots, l_b,$ and $m_0, m_1, \ldots, m_b$ of nonnegative integers which satisfy (\ref{eqn:rcondition}).
Using Lemmas \ref{lem:pathheightgf} and \ref{lem:qkdet}, we find that the generating function for step 2 is equal to
\begin{equation}
\label{eqn:step2}
\left( \frac{U_{k-2}\left(\frac{1-x}{2\sqrt{x}}\right)}{\sqrt{x} U_{k-1}\left(\frac{1-x}{2 \sqrt{x}}\right)} \right)^{l_0 - 1}
\end{equation}
and the generating functions for steps 3 and 4 are both equal to
\begin{equation}
\label{eqn:step34}
\frac{(-1)^{k-2}}{\sqrt{x} U_{k-1}\left(\frac{1-x}{2\sqrt{x}}\right)}.
\end{equation}
Taking the product of $x$, the quantities in (\ref{eqn:step1}) and (\ref{eqn:step2}), and the square of the quantity in (\ref{eqn:step34}), we obtain (\ref{eqn:12krgf}), as desired.
$\Box$

\bigskip

The following special case of Theorem \ref{thm:E2} is of particular interest.

\begin{corollary}
\label{cor:12k1gf}
For all $k \ge 2$ we have
\begin{equation}
\label{eqn:G12k}
\sum_{\pi} x^{|\pi|} = \frac{x}{\left( U_{k-1}\left(\frac{1-x}{2\sqrt{x}}\right) \right)^2}.
\end{equation}
Here the sum on the left is over all permutations in $\fS(1243, 2143)$ which contain exactly one subsequence of type $12\ldots k$.
\end{corollary}
\begin{proof}
Set $r = 1$ in Theorem \ref{thm:E2}.
\end{proof}

\section{Permutations in $\fS(1243, 2143)$ Which Avoid Another Pattern}
\label{sec:avoidenum}

In Section \ref{sec:Chebyshev} we used Theorem \ref{thm:pik} to obtain generating functions for $\fS_n(1243, 2143, \sigma 3 4 \ldots k)$, where $\sigma \in \fS_2$.
In this section we use Theorem \ref{thm:pik} to obtain generating functions for $\fS_n(1243, 2134, \sigma)$ for other interesting choices of $\sigma$;  in some cases we give explicit enumerations.
We begin with $\sigma = 231$.

\begin{proposition}
\label{prop:231}
For all $n \ge 2$,
\begin{equation}
\label{eqn:231}
|\fS_n(1243, 2143, 231)| = (n+2) 2^{n-3}.
\end{equation}
Moreover,
\begin{equation}
\label{eqn:231gf}
\sum_{n=0}^\infty |\fS_n(1243, 2143, 231)| x^n = 1 + x \frac{(x-1)^2}{(2x-1)^2}.
\end{equation}
\end{proposition}
\begin{proof}
We first prove (\ref{eqn:231gf}).
Let $S(x)$ denote the generating function on the left side of (\ref{eqn:231gf}).
Observe that if $\pi \in \fS_n(1243, 2143, 231)$ then $\pi^{-1}(n) = 1$, $\pi^{-1}(n) = 2$, or $\pi^{-1}(n) = n$.
Therefore we may partition $\fS(1243, 2143, 321)$ into the following four sets:
\begin{displaymath}
T_1 = \{\pi \in \fS(1243, 2143, 231)\ |\ |\pi| \le 3\};
\end{displaymath}
\begin{displaymath}
T_2 = \{\pi \in \fS(1243, 2143, 231)\ |\ |\pi| \ge 4 \hbox{ and } \pi(1) = |\pi|\};
\end{displaymath}
\begin{displaymath}
T_3 = \{\pi \in \fS(1243, 2143, 231)\ |\ |\pi| \ge 4 \hbox{ and } \pi(2) = |\pi|\};
\end{displaymath}
\begin{displaymath}
T_4 = \{\pi \in \fS(1243, 2143, 231)\ |\ |\pi| \ge 4 \hbox{ and } \pi(|\pi|) = |\pi|\}.
\end{displaymath}
The generating function for $T_1$ is $1 + x + 2x^2 + 5x^3$.
The map from $\fS_{n-1}(1243, 2143, 231)$ which carries $\pi$ to $n, \pi$ is a bijection between $\fS_{n-1}(1243, 2143, 231)$ and the set of permutations in $T_2$ which have length $n$, so the generating function for $T_2$ is $x (S(x) - 1 - x - x^2)$.
Now observe that if $n \ge 4$ and $\pi \in \fS_n(1243, 2143, 231)$ has $\pi(2) = n$ then $\pi(1) = 1$.
It follows that the map from $\fS_{n-2}(132, 231)$ to $\fS_n(1243, 2143, 231)$ which carries $\pi$ to $1,n,\hat{\pi}$, where $\hat{\pi}$ is the sequence obtained from $\pi$ by adding one to every entry, is a bijection between $\fS_{n-2}(132, 231)$ and the set of permutations in $T_3$ which have length $n$.
By \cite[Proposition 9]{SimionSchmidt}, the generating function for $\fS_n(132, 231)$ is $1 + {\displaystyle \frac{x}{1-2x}}$, so the generating function for $T_3$ is ${\displaystyle x^2 \left(1 + \frac{x}{1-2x} - 1 - x\right)}$.
Finally, the map from $\fS_{n-1}(1243, 2143, 231)$ to $\fS_n(1243, 2143, 231)$ which carries $\pi$ to $\pi,n$ is a bijection between $\fS_{n-1}(1243, 2143, 231)$ and the set of permutations in $T_4$ which have length $n$.
Therefore the generating function for $T_4$ is $x (S(x) - 1 - x - 2x^2)$.
Combine these observations to find that
$$S(x) = 1 + x + 2x^2 + 5x^3 + x (S(x) - 1 - x - 2x^2) + x^2 \left(1 + \frac{x}{1-2x} - 1 - x\right) + x (S(x) - 1 - x - 2x^2).$$
Solve this equation for $S(x)$ to obtain (\ref{eqn:231gf}).

Line (\ref{eqn:231}) is immediate from (\ref{eqn:231gf}).
\end{proof}

Combining Proposition \ref{prop:231} with Theorem \ref{thm:pik} and arguing as in the proof of Theorem \ref{thm:E1}, we find the generating function for $\fS_n(1243, 2143, 2314\ldots k)$ for any $k \ge 4$.

\begin{proposition}
Set $r_2(x) = (x-1)^2$ and 
$$r_k(x) = (1-2x)^2 (\sqrt{x})^{k-3} U_{k-3}\left(\frac{1-x}{2\sqrt{x}}\right) - (1-x)^2  (\sqrt{x})^{k-2} U_{k-4}\left(\frac{1-x}{2\sqrt{x}}\right) \hspace{30pt} (k \ge 3).$$
Then for all $k \ge 3$ we have
\begin{displaymath}
\sum_{n=0}^\infty |\fS_n(1243, 2143, 2314\ldots k)| x^n = 1 + x \frac{r_{k-1}(x)}{r_k(x)}.
\end{displaymath}
\end{proposition}

Next we consider $\sigma = 321$.
We begin with a lemma concerning which products $\pi_1 * \pi_2$ avoid 321.

\begin{lemma}
\label{lem:321}
Fix $i$ and $n$ such that $1 \le i \le n-2$, fix $\pi_1 \in \fS_i(1243, 2143, 321)$, and fix $\pi_2 \in \fS_{n-i}(1243, 2143, 321)$.
Then $\pi_1 * \pi_2 \in \fS_n(1243, 2143, 321)$ if and only if all of the following hold.
\renewcommand\labelenumi{{\upshape (\roman{enumi}) }}
\begin{enumerate}
\item
At least one of $\pi_1$ and $\pi_2$ begins with 1.
\item
The entries $2, 3, \ldots, i$ are in increasing order in $\pi_1$.
\item
The entries $\pi_2(2), \pi_2(3), \ldots, \pi_2(n-i)$ are in increasing order in $\pi_2$.
\end{enumerate}
\end{lemma}
\begin{proof}
($\Longrightarrow$)
Suppose $\pi_1 * \pi_2 \in \fS_n(1243, 2143, 321)$.
If (i) does not hold then the entries $(\pi_1 * \pi_2)(1), n-i, 1$ form a subsequence of type 321 in $\pi_1 * \pi_2$.
If (ii) does not hold then the smallest entry to the right of $n$, together with two of the entries among $2, 3, \ldots, i$ which are in decreasing order, form a subsequence of type 321 in $\pi_1 * \pi_2$.
If (iii) does not hold then the $n$, together with two of the entries among $\pi_2(2), \pi_2(3), \ldots, \pi_2(n-i)$ which are in decreasing order, form a subsequence of type 321 in $\pi_1 * \pi_2$.
Since we have a contradiction in each case, we find (i),(ii), and (iii) hold.

($\Longleftarrow$)
By Proposition \ref{prop:perminirecurrence}(i), the permutation $\pi_1 * \pi_2$ avoids 1243 and 2143.
Now observe that by (i)--(iii), one of the following holds.
\begin{enumerate}
\item
There exists $k$, $1 \le k \le i-1$, such that
$$\pi_1 * \pi_2 = n-i+1, n-i+2, \ldots, n-i+k, 1, n-i+k+1, \ldots, n-1,n,2,3,\ldots,n-i.$$
\item
There exists $k$, $1 \le k \le n-i$, such that
$$\pi_1 * \pi_2 = k, n-i+1, n-i+2, \ldots, n-1, n, 1, 2, \ldots, k-1, k+1, \ldots, n-i.$$
\end{enumerate}
It is routine to verify that in either case, $\pi_1 * \pi_2$ avoids 321.
Therefore $\pi_1 * \pi_2 \in \fS_n(1243, 2143, 321)$, as desired.
\end{proof}

We can now enumerate $\fS_n(1243, 2143, 321)$.

\begin{proposition}
\label{prop:321}
For all $n \ge 1$ we have
\begin{equation}
\label{eqn:321}
|\fS_n(1243, 2143, 321)| = {{n-1} \choose {0}} + {{n-1} \choose {1}} + 2 {{n-1} \choose {2}} + 2 {{n-1} \choose {3}}.
\end{equation}
Moreover,
\begin{equation}
\label{eqn:321gf}
\sum_{n=0}^\infty |\fS_n(1243, 2143, 321)| x^n = 1 + x \frac{1 - 2x + 3x^2}{(1-x)^4}.
\end{equation}
\end{proposition}
\begin{proof}
We first prove (\ref{eqn:321gf}).
Let $S(x)$ denote the generating function on the left side of (\ref{eqn:321gf}).
By Lemma \ref{lem:321}, we may partition $\fS(1243, 2143, 321)$ into the following four sets:
$$T_1 = \{\pi \in \fS(1243, 2143, 321)\ |\ \pi = \emptyset \hbox{ or } \pi(1) = n\};$$
$$T_2 = \{\pi \in \fS(1243, 2143, 321)\ |\ |\pi| \ge 2 \hbox{ and } \pi(|\pi|) = |\pi|\};$$
$$T_3 = \{\pi \in \fS(1243, 2143, 321)\ |\ \pi = \pi_1 * \pi_2, |\pi_2| \ge 2, \pi_1(1) = 1\};$$
$$T_4 = \{\pi \in \fS(1243, 2143, 321)\ |\ \pi = \pi_1 * \pi_2, |\pi_2| \ge 2, \pi_1(1) \neq 1, \pi_2(1) = 1\}.$$
Observe that $T_1$ consists of exactly those permutations of the form $n,1,2,3\ldots,n-1$, so the generating function for $T_1$ is ${\displaystyle \frac{1}{1-x}}$.
Observe that the map from $\fS_{n-1}(1243, 2143, 321)$ to $\fS_n(1243, 2143, 321)$ which carries $\pi$ to $\pi,n$ is a bijection from $\fS_{n-1}(1243, 2143, 321)$ to the set of permutations in $T_2$ which have length $n$, so the generating function for $T_2$ is $x (S(x) - 1)$.
Observe that if $\pi = \pi_1 * \pi_2$ and $\pi \in T_3$ then $\pi_1$ has the form $1,2,\ldots,|\pi_1|$ and $\pi_2$ has the form $k,1,2,\ldots, k-1, k+1, \ldots, |\pi_2|$ for some $k$, $1 \le k \le |\pi_2|$.
It follows that the generating function for $T_3$ is ${\displaystyle \left( \frac{x}{1-x}\right) \left( \frac{x}{(1-x)^2} - x\right)}$.
Finally, observe that if $\pi = \pi_1 * \pi_2$ and $\pi \in T_4$ then $\pi_2$ has the form $1,2,\ldots,|\pi_2|$ and there exists $k$, $2 \le k \le |\pi_1|$, such that $\pi_1$ has the form $2,3, \ldots, k, 1, k+1, \ldots, |\pi_1|$.
It follows that the generating function for $T_4$ is ${\displaystyle \left( \frac{x^2}{(1-x)^2}\right) \left( \frac{x^2}{1-x}\right)}$.
Combine these observations to find that
$$S(x) = \frac{1}{1-x} + x (S(x) - 1) + \frac{x^2}{(1-x)^3} - \frac{x^2}{1-x} + \frac{x^4}{(1-x)^3}.$$
Solve this equation for $S(x)$ to obtain (\ref{eqn:321gf}).

Line (\ref{eqn:321}) is immediate from (\ref{eqn:321gf}).
\end{proof}

Combining Proposition \ref{prop:321} with Theorem \ref{thm:pik} and arguing as in the proof of Theorem \ref{thm:E1}, we find the generating function for $\fS_n(1243, 2143, 3214\ldots k)$ for any $k \ge 4$.

\begin{proposition}
For all $k \ge 2$, set
$$r_k(x) = -(1+2x-x^2)(\sqrt{x})^{k+2} U_k\left( \frac{1-x}{2\sqrt{x}}\right) + (x^4-4x^3+2x^2+1)(\sqrt{x})^{k-1} U_{k-1}\left(\frac{1-x}{2\sqrt{x}}\right).$$
Then for all $k \ge 3$ we have
$$\sum_{n=0}^\infty |\fS_n(1243, 2143, 3214\ldots k)| x^n = 1 + x \frac{r_{k-1}(x)}{r_k(x)}.$$
\end{proposition}

At the end of Section \ref{sec:permproducts} we remarked that the results of that section can be generalized to $\fS(T_k)$, where $T_k$ is the set of all permutations in $\fS_k$ which end with $k,k-1$.
Our results in this section build almost exclusively on the results of Section \ref{sec:permproducts}, so they can also be generalized to $\fS(T_k)$.

\section{Permutations in $\fS(1243, 2143)$ Which Contain Another Pattern}

In Section \ref{sec:avoidenum} we used Theorem \ref{thm:pik} to obtain generating functions for $\fS_n(1243, 2143, \sigma)$ for various $\sigma$;  these results are elaborations on Theorem \ref{thm:E1}.
In this section we prove an analogue of Theorem \ref{thm:pik} for permutations in $\fS_n(1243, 2143)$ which contain a given pattern $\sigma$ exactly once.
We then use this result to prove several results which elaborate on Theorem \ref{thm:E2}.
We begin by setting some notation.

\begin{definition}
Fix a nonempty permutation $\tau \in \fS(1243, 2143)$ and set $k = |\tau| + 1$.
We write $G_\tau(x)$ to denote the generating function for those permutations in $\fS(1243, 2143)$ which contain exactly one subsequence of type $\tau$.
We write $H_\tau(x)$ to denote the generating function for those permutations in $\fS(1243, 2143)$ which avoid $\tau$.
We write $J_\tau(x)$ to denote the generating function for those permutations in $\fS(1243, 2143)$ which avoid $\tau,k$ and contain exactly one subsequence of type $\tau$.
\end{definition}

Our analogue of Theorem \ref{thm:pik} is the following result, which gives $G_{\tau,k}(x)$ in terms of $J_\tau(x)$ and $H_\tau(x)$.

\begin{theorem}
\label{thm:pik1}
Fix a nonempty permutation $\tau \in \fS(1243, 2143)$ and set $k = |\tau| + 1$.
Then
\begin{equation}
\label{eqn:pik1}
G_{\tau,k}(x) = \frac{x J_\tau(x)}{(2 - x - H_\tau(x))^2}.
\end{equation}
\end{theorem}
\begin{proof}
Observe that we may partition the set of permutations in $\fS(1243, 2143)$ which contain exactly one subsequence of type $\tau,k$ into the following three subsets.
$$T_1 = \{\pi\ |\ \pi \hbox{ begins with } |\pi|\};$$
$$T_2 = \{\pi\ |\ \pi = \pi_1 * \pi_2;  \pi_1 \hbox{ contains exactly one } \tau \hbox{ and no } \tau,k;  \pi_2 \hbox{ avoids } \tau,k\};$$
$$T_3 = \{\pi\ |\ \pi = \pi_1 * \pi_2;  \pi_1 \hbox{ avoids } \tau;  \pi_2 \hbox{ contains exactly one } \tau,k\}.$$
The generating function for $T_1$ is $x G_{\tau,k}(x)$, the generating function for $T_2$ is $J_\tau(x) (H_{\tau,k}(x) - 1)$, and the generating function for $T_3$ is $G_{\tau,k}(x) (H_\tau(x) - 1)$.
Combine these observations to find that
$$G_{\tau,k}(x) = x G_{\tau,k}(x) + J_\tau(x)(H_{\tau,k}(x) - 1) + G_{\tau,k}(x) (H_\tau(x) - 1).$$
Solve this last equation for $G_{\tau,k}(x)$ and use Theorem \ref{thm:pik} to eliminate $H_{\tau,k}(x)$.
Simplify the result to obtain (\ref{eqn:pik1}).
\end{proof}

Theorem \ref{thm:pik1} is particularly useful when $J_\tau(x) = G_\tau(x)$;  next we describe for which $\tau$ this occurs.

\begin{proposition}
\label{prop:G=J}
Fix a permutation $\tau \in \fS(1243, 2143)$ such that $|\tau| \ge 1$.
Then $J_\tau(x) = G_\tau(x)$ if and only if the last entry of $\tau$ is $|\tau|$.
\end{proposition}
\begin{proof}
($\Longrightarrow$)
Suppose $J_\tau(x) = G_\tau(x)$ and observe that this implies that if $\pi$ contains exactly one subsequence of type $\tau$ then $\pi$ avoids $\tau,|\tau|+1$.
However, if the last entry of $\tau$ is not $|\tau|$ then $\tau,|\tau|+1$ itself contains exactly one subsequence of type $\tau$ and does not avoid $\tau,|\tau|+1$.
This contradiction implies that the last entry of $\tau$ is $|\tau|$.

($\Longleftarrow$)
Suppose the last entry of $\tau$ is $|\tau|$.
If $\pi$ contains $\tau$ and $\tau, |\tau|+1$ then $\pi$ contains at least two subsequences of type $\tau$.
Therefore if $\pi$ contains exactly one subsequence of type $\tau$ then it avoids $\tau, |\tau|+1$;  it follows that $J_\tau(x) = G_\tau(x)$.
\end{proof}

Using Theorem \ref{thm:pik1} and Proposition \ref{prop:G=J}, one can compute $G_{\sigma}(x)$ for several infinite families of permutations.
Nevertheless, we content ourselves here with computing $G_{12\ldots k}(x)$ and $G_{213\ldots k}(x)$.
We begin by computing $G_{12\ldots k}(x)$, which amounts to giving another proof of Corollary \ref{cor:12k1gf}.

\bigskip
\noindent
{\it Another Proof of Corollary \ref{cor:12k1gf} }
We argue by induction on $k$.
To prove the result for $k = 2$, set $\tau = 1$ in (\ref{eqn:pik1}) and observe that $J_1(x) = x$ and $H_1(x) = 1$.
Now suppose the result holds for $k -1$.
Using (\ref{eqn:pik1}), Proposition \ref{prop:G=J}, (\ref{eqn:E112k}), and (\ref{eqn:qnUnrecurrence}) we find
\begin{eqnarray*}
G_{12\ldots k}(x) &=& \frac{x G_{12\ldots k-1}(x)}{\left( 1 - x - x \frac{U_{k-3}\left(\frac{1-x}{2\sqrt{x}}\right)}{\sqrt{x} U_{k-2}\left(\frac{1-x}{2\sqrt{x}}\right)}\right)^2} \\[2ex]
&=& \frac{x}{\left(\frac{1-x}{\sqrt{x}} U_{k-2}\left(\frac{1-x}{2\sqrt{x}}\right) - U_{k-3}\left(\frac{1-x}{2\sqrt{x}}\right)\right)^2} \\[3ex]
&=& \frac{x}{\left(U_{k-1}\left(\frac{1-x}{2\sqrt{x}}\right)\right)^2},
\end{eqnarray*}
as desired.
$\Box$

\bigskip

We now compute $G_{213\ldots k}(x)$.

\begin{proposition}
For all $k \ge 3$ we have
\begin{equation}
\label{eqn:G213k}
G_{213\ldots k}(x) = \frac{x (1+x) (1-x)^2}{\left(U_{k-1}\left(\frac{1-x}{2\sqrt{x}}\right)\right)^2}.
\end{equation}
\end{proposition}
\begin{proof}
First observe that the only permutations in $\fS(1243, 2143)$ which contain exactly one subsequence of type 21 and do not contain a subsequence of type 213 are 21 and 132, so $J_{21}(x) = x^2 + x^3$.
Now the result follows from (\ref{eqn:pik1}), Proposition \ref{prop:G=J}, and (\ref{eqn:qnUnrecurrence}) by induction on $k$.
\end{proof}

We conclude this section by showing that in a certain sense, there is no analogue of $\varphi$ which behaves as nicely with respect to $213\ldots k$ as $\varphi$ does with respect to $12\ldots k$.

\begin{corollary}
\label{cor:no213kbijection}
There is no bijection $\omega : \cS_n \longrightarrow \fS_{n+1}(1243, 2143)$ such that for every $\pi \in \cS_n$, the quantity $\tau_k(\pi)$ is equal to the number of subsequences of type $213\ldots k$ in $\omega(\pi)$.
\end{corollary}
\begin{proof}
If such a bijection existed then we would have $G_{12\ldots k}(x) = G_{213\ldots k}(x)$ for all $k \ge 3$, which contradicts (\ref{eqn:G12k}) and (\ref{eqn:G213k}).
\end{proof}

We remark that one can also prove Corollary \ref{cor:no213kbijection} by observing that if such a bijection existed, then there would be the same number of permutations in $\fS_4(1243, 2143)$ which contain exactly one subsequence of type 123 as there are containing exactly one subsequence of type 213.
However, there are 6 of the former (2314, 1423, 2341, 1342, 4123, 3124) and only 5 of the latter (3241, 2413, 1324, 3142, 4213).
As this example illustrates, subsequences of type $213\ldots k$ are more plentiful in $\fS_n(1243, 2143)$ than subsequences of type $12\ldots k$.

\section{Open Problems}

We conclude with some open problems suggested by the results in this paper.

\begin{enumerate}
\item
In view of (\ref{eqn:E112k}), (\ref{eqn:E121k}) and Corollary \ref{cor:varphiheight}, there exists a bijection $\rho : \cS_n \longrightarrow \fS_{n+1}(1243, 2143)$ such that for all $k \ge 2$, the Schr\"oder path $\pi$ has height $k-2$ or less if and only if $\rho(\pi)$ avoids $213\ldots k$.
Give a combinatorial description of such a bijection.
Then give a combinatorial description of the statistic $\alpha_k$, where $\alpha_k(\pi) = \tau_k(\rho^{-1}(\pi))$ for all $\pi \in \fS_n(1243, 2143)$.
Alternatively, give a bijection $\omega : \fS_n(1243, 2143) \longrightarrow \fS_n(1243, 2143)$ such that for all $k \ge 0$ and any permutation $\pi \in \fS_n(1243, 2143)$, we have that $\pi$ avoids $12\ldots k$ if and only if $\omega(\pi)$ avoids $213\ldots k$. 

\item
For any permutation $\pi \in \fS_n(1243, 2143)$ let $\cT_\pi(x)$ denote the generating function given by 
$$\cT_\pi(x) = \sum_{n=0}^\infty |\fS_n(1243, 2143, \pi)| x^n.$$
Then the case of Theorem \ref{thm:pik} in which $|R| = 1$ amounts to a recurrence relation giving $\cT_{\pi * 1}(x)$ in terms of $\cT_\pi(x)$.
Find similar relations for $\cT_{\pi_1 * \pi_2}(x)$ in terms of $\cT_{\pi_1}(x)$ and $\cT_{\pi_2}(x)$ and for $\cT_{n\pi}(x)$ in terms of $\cT_\pi(x)$.
Alternatively, find an analogue of \cite[Theorem 2.1]{MansourVainshtein2} for $\fS(1243, 2143)$.

\item
Give a combinatorial proof of the fact that
$$|\fS_n(1243, 2143, 231)| = (n+2) 2^{n-3} \hspace{30pt} (n \ge 2).$$

\item
Give a combinatorial proof of the fact that
$$|\fS_n(1243, 2143, 321)| = {{n-1} \choose {0}} + {{n-1} \choose {1}} + 2 {{n-1} \choose {2}} + 2 {{n-1} \choose {3}} \hspace{30pt} (n \ge 1).$$

\item
For $k \ge 3$ and $r \ge 2$, find the generating function for the number of permutations in $\fS_n(1243, 2143)$ which contain exactly $r$ subsequences of type $213\ldots k$.
\end{enumerate}

\bigskip

\begin{Large}
\noindent
{\bf Acknowledgment}
\end{Large}

\bigskip

The first author thanks Justin Pryzby, Gettysburg College class of 2005, for generating data which led to Theorem \ref{thm:contfrac}.

\end{document}